\documentclass[preprint,12pt]{elsarticle}

\usepackage{amssymb}
\usepackage{amsmath}
\usepackage{amsthm}
\usepackage{subcaption}

\newtheorem{theorem}{Theorem}[section]
\newtheorem{lemma}[theorem]{Lemma}

\theoremstyle{remark}
\newtheorem{remark}[theorem]{Remark}

\newcounter{roughpanel}

\newcommand{\resultpanel}[2]{%
    \begin{minipage}[t]{0.31\textwidth}
        \centering
        \includegraphics[width=\linewidth]{#1}
        \par\vspace{2pt}
        \refstepcounter{roughpanel}
        \small (\alph{roughpanel})~#2
    \end{minipage}%
}

\newcounter{widepanel}

\newcommand{\resultpanelwide}[2]{%
    \begin{minipage}[t]{0.47\textwidth}
        \centering
        \includegraphics[width=\linewidth]{#1}
        \par\vspace{2pt}
        \refstepcounter{widepanel}
        \small (\alph{widepanel})~#2
    \end{minipage}%
}

\journal{Journal of Computational Physics}

\begin{document}

\begin{frontmatter}



\title{A time-domain direct sampling method for inverse scattering problems by rough surfaces}


\author[label1]{Xiaoli Liu}
\ead{xiaoli_liu@buaa.edu.cn}
\author[label1]{Yun Liu}
\ead{liuyyun@buaa.edu.cn}
\author[label2]{Haiwen Zhang}
\ead{zhanghaiwen@amss.ac.cn}


\affiliation[label1]{organization={School of Mathematical Sciences},
            addressline={Beihang University}, 
            city={Beijing},
            postcode={100191}, 
            country={China}}
\affiliation[label2]{organization={SKLMS and Academy of Mathematics and Systems Science},
            addressline={Chinese Academy of Sciences}, 
            city={Beijing},
            postcode={100190}, 
            country={China}}

\begin{abstract}
    In this paper, we propose a time-domain direct sampling method for reconstructing a two-dimensional unbounded rough surface from transient near-field measurements. We first establish the well-posedness of the time-domain forward scattering problem by using the Laplace transform, the Dirichlet-to-Neumann operator, and variational methods. Based on a time-domain Helmholtz--Kirchhoff identity, we construct a back-propagated field and show that it can be interpreted as the scattered field generated by a virtual incident wave. To improve the focusing behavior of the imaging method, a correction term is further introduced to compensate for the missing lower-aperture contribution, leading to a direct imaging functional for locating the rough surface. The proposed method makes use of the broadband information contained in a single time-domain pulse without frequency sweeping. Numerical experiments are presented to demonstrate the effectiveness of the proposed method in reconstructing rough surfaces with different geometric features. Comparisons with frequency-domain direct sampling methods are also provided, together with investigations of the resolution, the influence of measurement configuration, and the robustness with respect to noise.
\end{abstract}




\end{frontmatter}



\section{Introduction}\label{sec:introduction}

Inverse scattering problems by unbounded rough surfaces have attracted considerable attention due to their important applications in various scientific and engineering fields, such as geophysical exploration, underwater acoustics, remote sensing, and nondestructive testing. In this paper, we consider the inverse acoustic scattering problem by a two-dimensional unbounded rough surface in the time domain. Our aim is to reconstruct the unknown rough surface from transient near-field measurements generated by incident acoustic pulses.

Various numerical methods have been developed for inverse scattering problems by unbounded rough surfaces, mainly in the frequency domain. For iterative optimization methods, Newton-type algorithms and related approaches have been studied to reconstruct locally rough surfaces and periodic interfaces \cite{BaoLin2011,LiuZhang2019Periodic}. 
These methods can provide accurate reconstructions of the unknown surfaces, but their performance may depend on a suitable initial guess, and repeated solutions of the
forward scattering problem often lead to relatively high computational costs. The transformed field expansion method has also been developed for reconstructing small and smooth
perturbations of planar surfaces \cite{BaoLi2013,BaoLi2014}. On the other hand, several
non-iterative sampling-type methods have been proposed, including the linear sampling method and the factorization method for locally rough or periodic surfaces
\cite{ArensKirsch2003,DingEtAl2017,LechleiterNguyen2013,LiYangZhang2022}. In particular, direct imaging methods have been developed for acoustic and elastic scattering by general unbounded rough surfaces \cite{HuEtAl2019,LiuZhangZhang2018,LiuZhangZhang2019}. These methods construct imaging functionals directly from the measured data and generally
do not require iterative solutions of the forward problem. More recently, direct imaging and reverse time migration methods have also been investigated for locally rough interfaces \cite{LiEtAl2024,LiWuYang2025}. These studies provide important foundations for rough surface reconstruction, but are mainly concerned with time-harmonic data.

Inverse scattering problems using time-domain measurements have also been extensively investigated. Compared with time-harmonic data, transient measurements retain the temporal information of wave propagation and contain broadband information from a single incident pulse.
The well-posedness of time-domain forward scattering problems has been studied using various analytical approaches. In particular, Bamberger and Ha-Duong \cite{BambergerHaDuong1986} developed a space-time variational formulation based on retarded potentials, while Sayas \cite{Sayas2016} provided a systematic framework for time-domain boundary integral equations.
For scattering in unbounded structures, Gao et al. \cite{GaoLiZhang2017} established well-posedness and stability estimates for transient acoustic-elastic wave scattering by using transparent boundary conditions, the Laplace transform, and variational methods.
Several sampling-type methods have been developed directly in the time domain. In particular, the time-domain linear sampling method was introduced for inverse obstacle scattering problems and subsequently extended to other boundary conditions \cite{ChenEtAl2010,HaddarLechleiterMarmorat2014}. The factorization method has also been studied in the time domain for different scattering configurations \cite{CakoniHaddarLechleiter2019,HaddarLiu2020}. For rough surface reconstruction, a time-domain probe method was proposed in \cite{BurkardPotthast2009}. More closely related to the present work, direct sampling and imaging methods have been developed to reconstruct
acoustic scatterers from transient measurements \cite{GuoEtAl2016,GuoLiWang2024}. Further developments include time-domain sampling methods with experimental measurements and passive data \cite{LiuSongEtAl2023,LiuTianZhang2026}. These studies provide a basis for developing time-domain imaging methods for different scattering configurations.

Motivated by the above studies, it is of interest to develop a direct sampling method for reconstructing general unbounded rough surfaces using time-domain measurements. However, the unbounded geometry introduces additional difficulties in extending existing imaging
methods for bounded scatterers. In particular, the rough surface cannot be enclosed by a finite measurement curve, and the derivation of the corresponding time-domain integral identities requires a careful treatment of wave propagation in unbounded domains. 
A suitable well-posedness framework for the time-domain forward scattering problem is also needed to support the subsequent imaging analysis.
Another important issue is how to construct a back-propagated field from transient near-field data and develop a direct imaging functional with a clear mathematical interpretation of its imaging mechanism and focusing behavior. Moreover, since the measurements are taken only above the rough surface, a correction term can be introduced to compensate for the missing lower-aperture contribution and improve the focusing behavior of the imaging functional.

In this paper, we propose a time-domain direct sampling method for reconstructing an unbounded rough surface from transient near-field measurements. We first establish the well-posedness of the associated forward scattering problem using the Laplace transform, a Dirichlet-to-Neumann operator, and variational arguments. We then derive a time-domain Helmholtz--Kirchhoff identity suited to the unbounded geometry and construct a back-propagated field from the measured Cauchy data. The field is related to an auxiliary scattering problem with a virtual incident wave, which provides a mathematical interpretation of the imaging process. To improve the focusing behavior, we introduce a correction term associated with the missing lower-aperture contribution and construct a direct imaging functional. We further investigate the focusing properties of the corresponding idealized full-aperture kernel and discuss the imaging mechanism. Numerical experiments with different rough surface profiles demonstrate the effectiveness of the method. Comparisons with single- and multi-frequency direct sampling methods are presented, and the effects of pulse center-frequency, measurement configuration, and noise are investigated.

The remaining part of this paper is organized as follows. Section~2 introduces the mathematical formulation and preliminary results. Section~3 establishes the well-posedness of the time-domain forward scattering problem. In Section~4, we derive the time-domain Helmholtz--Kirchhoff identity, construct the imaging functional, and analyze its imaging mechanism. Section~5 presents numerical experiments to validate the method and investigate its performance. Conclusions and possible future directions are given in Section~6.

\section{Mathematical formulation}
We consider the two-dimensional time-domain acoustic scattering problem by an unbounded impenetrable rough surface subject to a Dirichlet boundary condition.
Denote by {$\Gamma := \{(x_1,x_2)\in \mathbb{R}^2 : x_2 = f_\Gamma(x_1)\}$} the rough surface, $D^+ := \{(x_1,x_2)\in \mathbb{R}^2 : x_2 > f_\Gamma(x_1)\}$ the region above the surface, where {$f_\Gamma$ is a bounded Lipschitz function}.
Let $p^s(x,t)$ denote the scattered acoustic pressure field generated by a causal incident field $p^i(x,t)$, which means $p^i(x,t) = 0$ for $t<0$.
The scattering problem is formulated as
\begin{equation}\label{eq:problem_original}
    \begin{cases}
        c^{-2} \partial_{tt} p^s(x,t) - \Delta_x p^s(x,t) = 0, & (x,t) \in D^+ \times \mathbb{R}, \\
        p^s(x,t) = -p^i(x,t), & (x,t) \in \Gamma \times \mathbb{R}, \\
        p^s(x,t) = \partial_t p^s(x,t) = 0, & (x,t) \in D^+ \times \mathbb{R}\backslash\mathbb{R}^+,
    \end{cases}
\end{equation}
where $c$ is the wave speed. The causal fundamental solution of the two-dimensional wave equation is 
\begin{equation*}
    \Phi(x,t;y):=\frac{{H}(t-|x-y|/c)}{2\pi\sqrt{t^2-c^{-2}|x-y|^2}},\qquad (x,t)\in\left\{\mathbb{R}^2\backslash\{y\}\right\}\times\mathbb{R},
\end{equation*}
where $H$ is the Heaviside function and $\Phi$ satisfies
\begin{equation*}
    c^{-2}\partial_{tt}\Phi(x,t;y)-\Delta \Phi(x,t;y) = \delta(x-y)\delta(t),\qquad (x,t)\in\left\{\mathbb{R}^2\right\}\times\mathbb{R}.  
\end{equation*}

We denote by $\mathcal{D} (\mathbb{R};X) = C_0^\infty(\mathbb{R};X)$ the set of smooth and compactly supported $X$-valued functions for a Hilbert space $X$.
Further, $\mathcal{D}' (\mathbb{R} ; X)$ are $X$-valued distributions on the real line, and $\mathcal{S}' (\mathbb{R} ; X)$ are the corresponding tempered distributions.
For a causal temporal signal $\chi(t)\in\mathcal{S}' (\mathbb{R} ; X)$, we consider an incident wave field 
\begin{equation}\label{inc}
    p^i(x,t;y):=\left[\chi*\Phi(x,\cdot;y) \right](t),\qquad (x,t)\in\left\{\mathbb{R}^2\backslash\{y\}\right\}\times\mathbb{R},
\end{equation}
where {$*$} denotes convolution in time.

\begin{figure}[htbp]
    \centering
    \includegraphics[width=0.6\linewidth]{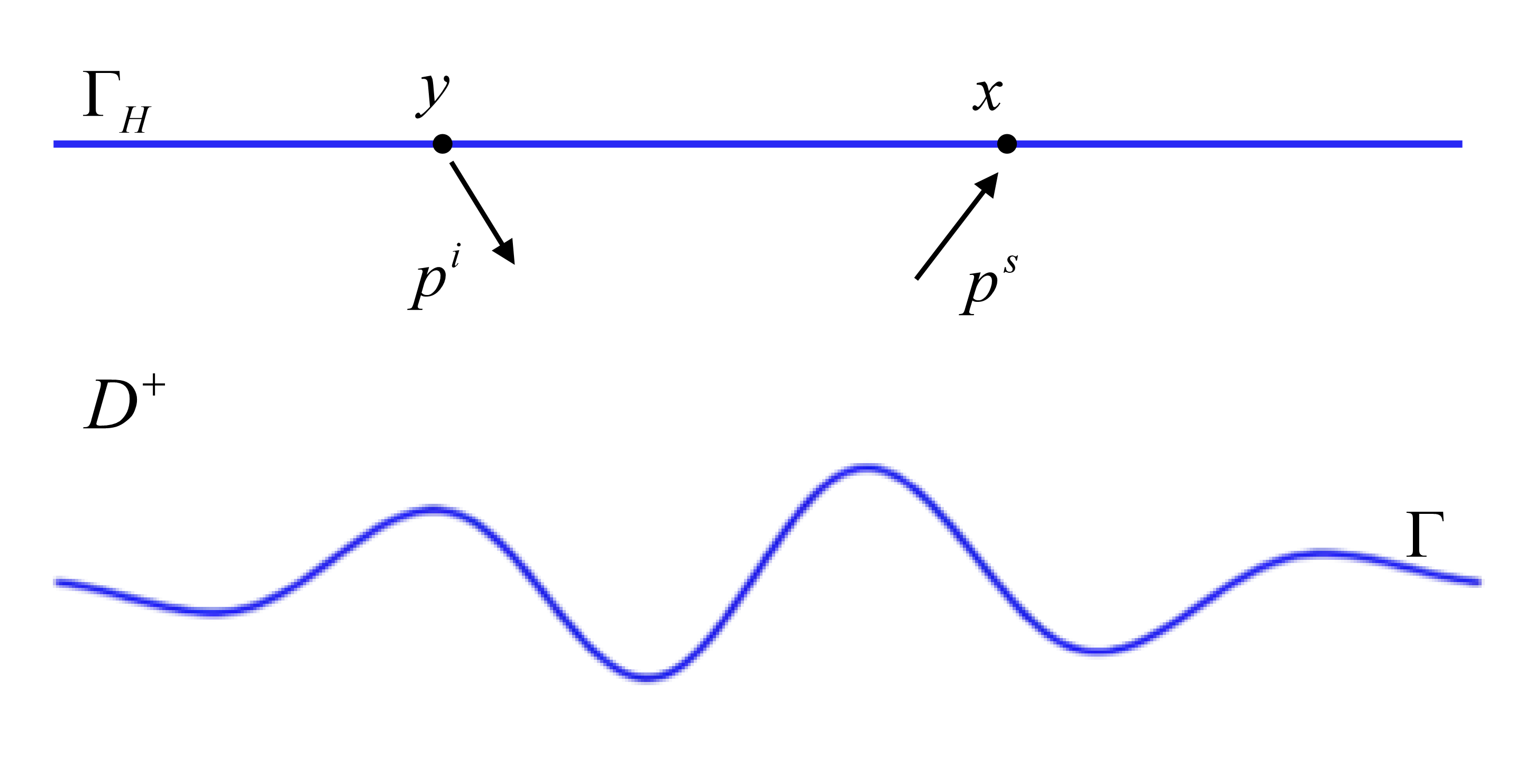}
    \caption{Geometry setting for the scattering problem}
    \label{fig:geometry_setting}
\end{figure}

The \emph{forward scattering problem} is to find a solution $p^s$ satisfying \eqref{eq:problem_original} corresponding to an incident wave $p^i$ in \eqref{inc}. For the inverse problem, let
\begin{equation*}
\Gamma_H:=\{(x_1,H) : x_1\in\mathbb R\},
\qquad
H>\sup_{x_1\in\mathbb R}f_\Gamma(x_1),
\end{equation*}
be a horizontal measurement line located above the rough surface. The \emph{inverse scattering problem} is to reconstruct the unknown rough surface $\Gamma$ from the time-domain Cauchy data 
\begin{equation*}
\left\{\left(p^s(x,t;y),\partial_\nu p^s(x,t;y)\right) : (x,t,y)\in \Gamma_H\times\mathbb{R}\times\Gamma_H \right\},
\end{equation*}
where $\partial_\nu$ is the normal derivative on $\Gamma_H$ pointing to the upper half-space. Here $y\in\Gamma_H$ denotes the source location, while $x\in\Gamma_H$ denotes the observation point. The measurement configuration is illustrated in Figure \ref{fig:geometry_setting}.

\section{Well-posedness of the forward scattering problem}

In this section, we consider the well-posedness of the forward scattering problem \eqref{eq:problem_original}, aiming to establish a theoretical foundation for the subsequent inversion algorithm. The well-posedness of time-domain rough surface scattering problems has been studied in the literature; see, e.g., \cite{GaoLiZhang2017}. In particular, the authors of \cite{GaoLiZhang2017} employed a transparent boundary condition, the Laplace transform, and variational methods to establish well-posedness and stability estimates for transient wave scattering in an unbounded structure.

Following the analytical framework of \cite{GaoLiZhang2017}, we employ the Laplace transform and variational methods to derive the corresponding well-posedness results for the time-domain acoustic rough surface problem considered in this work. The original problem \eqref{eq:problem_original} can be equivalently reformulated as the following reduced problem:
\begin{equation}\label{eq:problem_reduced}
    \begin{cases}
        c^{-2}\partial_{tt} p^s(x,t) - \Delta_x p^s(x,t) = 0, \qquad & (x,t) \in D_h \times \mathbb{R}, \\
        p^s(x,t) = -p^i(x,t), \qquad & (x,t) \in \Gamma \times \mathbb{R}, \\
        p^s(x,t) = \partial_t p^s(x,t) = 0, \qquad & (x,t) \in D_h \times \mathbb{R}\backslash\mathbb{R}^+, \\
        \partial_\nu p^s(x,t) = \mathcal{J}p^s(x,t), \qquad & (x,t) \in \Gamma_h \times \mathbb{R},
    \end{cases}
\end{equation}
where $\Gamma_h:=\{(x_1,x_2)\in\mathbb{R}^2:x_2=h\}$ is a straight line above the rough surface, with
\[
\sup_{x_1\in\mathbb R} f_\Gamma(x_1)<h<H,
\]
and $D_h:=\{(x_1,x_2)\in\mathbb{R}^2:f_\Gamma(x_1)<x_2<h\}$ denotes the strip-like domain between $\Gamma$ and $\Gamma_h$.
The operator $\mathcal{J}:=\mathcal{L}^{-1}\circ\mathcal{B}\circ\mathcal{L}$, where $\mathcal{L}$ and $\mathcal{B}$ denote the Laplace transform and the Dirichlet-to-Neumann (DtN) operator, respectively.

For $\sigma\in\mathbb R$, we set
\[
\mathcal{L}'_\sigma(\mathbb{R};X)
=
\left\{
f\in\mathcal{D}'(\mathbb{R};X):
e^{-\sigma t}f(t)\in\mathcal{S}'(\mathbb{R};X)
\right\}.
\]
For a causal function $f$ admitting a Laplace transform, we define
\[
\breve f(s)=\mathcal{L}[f](s)
:=
\int_{-\infty}^{+\infty}f(t)e^{\mathrm{i}st}\,\mathrm{d}t,
\qquad s=\omega+\mathrm{i}\sigma,
\]
where $\sigma>0$ and $\omega\in\mathbb R$.

For $u\in H^{1/2}(\Gamma_h)$, the DtN operator $\mathcal{B}(s)$ is defined by
\[
(\mathcal{B}(s)u)(x_1)
=
-\frac{1}{2\pi}
\int_{\mathbb R}
\beta(\xi,s)\hat u(\xi)e^{-\mathrm{i}\xi x_1}
\,\mathrm{d}\xi,
\]
where
\[
\beta(\xi,s)^2
=
\xi^2-\frac{s^2}{c^2},
\qquad
\operatorname{Re}\beta(\xi,s)>0,
\]
and
\[
\hat u(\xi)
=
\int_{\mathbb R}
u(x_1)e^{\mathrm{i}\xi x_1}\,\mathrm{d}x_1.
\]

To obtain the well-posedness of the reduced problem \eqref{eq:problem_reduced}, we first consider the related problem in the Laplace domain. Let $\breve{p^s}(x,s)$ and $\breve{p^i}(x,s)$ denote the Laplace transforms of $p^s(x,t)$ and $p^i(x,t)$ with respect to $t$, respectively.

Taking the Laplace transform of \eqref{eq:problem_reduced} and using the zero initial conditions, we obtain
\begin{equation}\label{eq:problem_Laplace-domain}
    \begin{cases}
        \Delta_x\breve{p^s}(x,s)+s^2c^{-2}\breve{p^s}(x,s)=0,
        & x\in D_h,\\
        \breve{p^s}(x,s)=-\breve{p^i}(x,s),
        & x\in\Gamma,\\
        \partial_\nu\breve{p^s}(x,s)
        =\mathcal{B}(s)\breve{p^s}(x,s),
        & x\in\Gamma_h.
    \end{cases}
\end{equation}


Before presenting the well-posedness result for problem
\eqref{eq:problem_Laplace-domain}, we recall some useful
properties of the DtN operator $\mathcal{B}(s)$.

By the standard trace theorem and
\cite[Lemmas 2.2, 2.4 and 2.5]{GaoLiZhang2017},
there exists a positive constant $C_{\mathrm{tr}}$
such that
\begin{equation*}
    \|v\|_{H^{1/2}(\Gamma_h)}
    \leq C_{\mathrm{tr}}\|v\|_{H^1(D_h)},
    \qquad \forall v\in H^1(D_h).
\end{equation*}
Moreover, for $s=\omega+\mathrm{i}\sigma$ with
$\sigma\geq\sigma_0>0$, the DtN operator
$\mathcal{B}(s):H^{1/2}(\Gamma_h)
\to H^{-1/2}(\Gamma_h)$ satisfies
\begin{equation}
    \|\mathcal{B}(s)v\|_{H^{-1/2}(\Gamma_h)}
    \leq C(\sigma_0)|s|
    \|v\|_{H^{1/2}(\Gamma_h)},
    \label{eq:DtN_continuity}
\end{equation}
and
\begin{equation}
    -\mathrm{Re}\,
    \left\langle
    \mathrm{i}s^{-1}\mathcal{B}(s)v,v
    \right\rangle_{\Gamma_h}
    \geq 0,
    \label{eq:DtN_sign}
\end{equation}
for all $v\in H^{1/2}(\Gamma_h)$,
where $C(\sigma_0)=\max\{c^{-1},\sigma_0^{-1}\}$.

With these properties, we establish the following
well-posedness result in the Laplace domain.

\begin{theorem}\label{thm:well-posedness_Laplace-domain}
    Let $s=\omega+\mathrm{i}\sigma$ with $\sigma>0$.
    For any $\breve{p^i}|_\Gamma\in H^{1/2}(\Gamma)$,
    problem \eqref{eq:problem_Laplace-domain} has a unique
    weak solution $\breve{p^s}\in H^1(D_h)$, which satisfies
    \begin{equation}\label{eq:stabilityEstimate_Laplace-domain}
        \|\nabla\breve{p^s}\|_{L^2(D_h)^2}
        +\|s\breve{p^s}\|_{L^2(D_h)}
        \leq C(s)\|\breve{p^i}\|_{H^{1/2}(\Gamma)},
    \end{equation}
    where $C(s)>0$ is a constant independent of
    the boundary data.
\end{theorem}

\begin{proof}
Since $\breve{p^s}|_\Gamma=-\breve{p^i}|_\Gamma
    \in H^{1/2}(\Gamma)$, by the trace extension theorem,
    there exists $g\in H^1(D_h)$ such that
     \begin{equation}\label{g}
        g|_\Gamma = \breve{p^s}|_\Gamma\qquad\text{and}\qquad \|g\|_{H^1(D_h)} \leq C_1 \|\breve{p^s}\|_{H^{1/2}(\Gamma)}.
    \end{equation}
    Define
    \[
        H^1_\Gamma(D_h)
        =\{v\in H^1(D_h):v|_\Gamma=0\}
    \]
    and set $u=\breve{p^s}-g$. Then
    $u\in H^1_\Gamma(D_h)$ satisfies the variational problem
    \begin{equation}\label{eq:problem_variational}
        a_s(u,v)=-a_s(g,v),
        \qquad \forall v\in H^1_\Gamma(D_h),
    \end{equation}
    where
    \begin{align*}
        a_s(u,v)
        ={}
        \int_{D_h}
        \left(
        \mathrm{i}s^{-1}\nabla u\cdot\nabla\overline v
        -\mathrm{i}sc^{-2}u\overline v
        \right)\mathrm{d}x-
        \left\langle
        \mathrm{i}s^{-1}\mathcal{B}(s)u,v
        \right\rangle_{\Gamma_h}.
    \end{align*}

    By the trace theorem and the continuity of
    $\mathcal{B}(s)$, the sesquilinear form $a_s$
    is bounded on $H^1_\Gamma(D_h)$, namely,
    \[
        |a_s(u,v)|
        \leq C(s)\|u\|_{H^1(D_h)}
                    \|v\|_{H^1(D_h)}.
    \]

    Moreover, using the sign property of the DtN
    operator, we obtain
    \begin{align*}
        \mathrm{Re}\,a_s(v,v)
        &\geq
        \frac{\sigma}{|s|^2}
        \|\nabla v\|_{L^2(D_h)^2}^2
        +\frac{\sigma}{c^2}
        \|v\|_{L^2(D_h)}^2\\
        &\geq
        \sigma\min\{|s|^{-2},c^{-2}\}
        \|v\|_{H^1(D_h)}^2.
    \end{align*}
    Thus, $a_s$ is coercive for each fixed
    $s$ with $\mathrm{Im}\,s>0$.

    It follows from the Lax--Milgram theorem that
    \eqref{eq:problem_variational} admits a unique
    solution $u\in H^1_\Gamma(D_h)$. Consequently,
    $\breve{p^s}=u+g$ is the unique weak solution
    of \eqref{eq:problem_Laplace-domain}.

    Finally, taking $v=u$ in
    \eqref{eq:problem_variational}, the boundedness
    and coercivity of $a_s$ yield
    \[
        \|u\|_{H^1(D_h)}
        \leq C(s)\|g\|_{H^1(D_h)}.
    \]
    Combining this estimate with \eqref{g} and
    the triangle inequality gives
    \[
        \|\breve{p^s}\|_{H^1(D_h)}
        \leq C(s)
        \|\breve{p^i}\|_{H^{1/2}(\Gamma)}.
    \]
    Since $s$ is fixed, the stability estimate
    \eqref{eq:stabilityEstimate_Laplace-domain}
    follows immediately.
\end{proof}

The above theorem establishes the well-posedness of the forward scattering problem in the Laplace domain. We now consider the corresponding problem in the time domain for the pulsed point source defined in \eqref{inc}.

\begin{theorem}\label{thm:well-posedness_reduced}
Let $y=(y_1,H)\in\Gamma_H$ and assume that $\sup f_\Gamma<h<H$. Suppose that the incident field is given by
\[
    p^i(x,t;y)
    =[\chi*\Phi(x,\cdot;y)](t),
\]
where $\chi\in H^1(\mathbb R)$ is causal. Then, for any $T>0$, the reduced problem
\eqref{eq:problem_reduced} admits a unique solution
\[
    p^s\in L^2(0,T;H^1(D_h))
    \cap H^1(0,T;L^2(D_h)).
\]
Moreover, the solution satisfies the stability estimate
\begin{equation}\label{eq:stabilityEstimate_reduced}
    \|p^s\|_{L^2(0,T;H^1(D_h))
    \cap H^1(0,T;L^2(D_h))}
    \leq C_T\|\chi\|_{H^1(\mathbb R)},
\end{equation}
where $C_T>0$ is independent of $\chi$.
\end{theorem}
\begin{proof}
Let $s=\omega+\mathrm{i}\sigma$ with $\sigma>0$, and denote by $\breve p=\breve p^i+\breve p^s$ the total field in the Laplace domain. Thus, $\breve p=0$ on $\Gamma$
and
\[
    \partial_\nu\breve p
    =\mathcal B(s)\breve p+F_s
    \quad\text{on }\Gamma_h,
\]
where
\[
    F_s:=\partial_\nu\breve p^i
    -\mathcal B(s)\breve p^i.
\]
Using the variational form introduced in the proof
of Theorem~\ref{thm:well-posedness_Laplace-domain},
we obtain
\[
    a_s(\breve p,v)
    =\left\langle
      \mathrm{i}s^{-1}F_s,v
     \right\rangle_{\Gamma_h},
    \qquad v\in H^1_\Gamma(D_h).
\]

By the coercivity of $a_s$, we have
\[
    \mathrm{Re}\,a_s(v,v)
    \geq \frac{\sigma}{|s|^2}
    \left(
        \|\nabla v\|_{L^2(D_h)^2}^2
        +\frac{|s|^2}{c^2}
         \|v\|_{L^2(D_h)}^2
    \right).
\]
Moreover, since $v=0$ on $\Gamma$, the
one-dimensional trace inequality yields
\[
    \|v\|_{L^2(\Gamma_h)}^2
    \leq 2\|v\|_{L^2(D_h)}
            \|\partial_{x_2}v\|_{L^2(D_h)}
    \leq \frac{2c}{|s|}
    \left(
        \|\nabla v\|_{L^2(D_h)^2}^2
        +\frac{|s|^2}{c^2}
         \|v\|_{L^2(D_h)}^2
    \right).
\]
Taking $v=\breve p$ in the variational equation
and using the above inequalities, we obtain
\[
\begin{aligned}
\left(
    \|\nabla\breve p\|_{L^2(D_h)^2}^2
    +\frac{|s|^2}{c^2}
     \|\breve p\|_{L^2(D_h)}^2
 \right)^{1/2}
\leq
    \frac{\sqrt{2c}|s|^{1/2}}{\sigma}
    \|F_s\|_{L^2(\Gamma_h)}.
\end{aligned}
\]

Fix $\sigma=\sigma_0>0$.
By Lemma~\ref{lem:point_source_bounds}
in \ref{app:point_source_estimates},
we have
\[
    \|F_s\|_{L^2(\Gamma_h)}
    +\left(
        \|\nabla\breve p^i\|_{L^2(D_h)^2}^2
        +\frac{|s|^2}{c^2}
         \|\breve p^i\|_{L^2(D_h)}^2
     \right)^{1/2}
    \leq C_{\sigma_0}(1+|s|)^{1/2}
    |\breve\chi(s)|.
\]
Since $\breve p^s=\breve p-\breve p^i$,
the triangle inequality gives
\begin{equation}\label{eq:point_source_frequency_estimate}
    \|\nabla\breve p^s\|_{L^2(D_h)^2}
    +\|s\breve p^s\|_{L^2(D_h)}
    \leq C_{\sigma_0}(1+|s|)
    |\breve\chi(s)|.
\end{equation}
As shown in \ref{app:point_source_estimates},
$\breve p^s$ is analytic in the upper half-plane
and satisfies the required polynomial growth estimate.
Hence, by {\cite[Lemma 43.1]{Treves1975}},
the inverse Laplace transform of $\breve p^s$
exists and is supported in $[0,\infty)$.
By Parseval's identity,
\begin{align*}
 &\int_0^\infty e^{-2\sigma_0t}
 \big(\|\nabla p^s(t)\|_{L^2(D_h)}^2+
 \|\partial_t p^s(t)\|_{L^2(D_h)}^2\big)\,\mathrm dt\\
 &\quad\le C_{\sigma_0}
 \int_{\mathbb R}(1+|s|)^2
 |\breve\chi(s)|^2\,\mathrm d\omega
 \le C_{\sigma_0}\|\chi\|_{H^1(\mathbb R)}^2.
\end{align*}
Here we used causality and
$\breve\chi(\omega+\mathrm{i}\sigma_0)
=\mathcal F(e^{-\sigma_0t}\chi)(\omega)$.
Finally, $p^s(0)=0$ implies
$\|p^s\|_{L^2(0,T;L^2(D_h))}
\le T\|\partial_t p^s\|_{L^2(0,T;L^2(D_h))}$.
Removing the exponential weight on $(0,T)$
proves \eqref{eq:stabilityEstimate_reduced}.
Uniqueness follows from the Laplace-domain
uniqueness.
\end{proof}

\section{Imaging algorithm}

In this section, we develop a time-domain direct sampling method for the inverse scattering problem by a two-dimensional rough surface. The construction follows the chain from measured transient Cauchy data, through a time-domain Helmholtz--Kirchhoff identity and back-propagation procedure, to a virtual scattering field and a corrected imaging functional. Unlike frequency-domain direct sampling methods, the present approach utilizes transient wave fields directly and therefore naturally incorporates broadband information from a single pulse signal without frequency sweeping.

\subsection{Time-domain Helmholtz--Kirchhoff identity}\label{TDHK}

We first establish a time-domain Helmholtz--Kirchhoff identity on the measurement line. This identity converts a Cauchy integral over $\Gamma_H$ into an aperture kernel, which is the analytical basis for the subsequent back-propagation construction.

Let
\[
\Gamma_H:=\{(x_1,H):x_1\in\mathbb R\}
\]
be the measurement line and
\[
U_H:=\{x\in\mathbb R^2:x_2>H\}.
\]
Define
\[
S_+^1
=
\{\hat x\in\mathbb R^2:|\hat x|=1,\ \hat x_2>0\},
\]
and introduce the half-aperture Helmholtz--Kirchhoff kernel
\[
K(t;y,z)
:=
\frac1{4\pi^2}
\int_{S_+^1}
\operatorname{p.v.}
\frac1{
t-c^{-1}\hat x\cdot(z-y)
}
\,\mathrm ds(\hat x).
\]

\begin{lemma}[Time-domain Helmholtz--Kirchhoff identity]
Assume that $y,z\in \mathbb R^2\setminus \overline{U_H}$.
Let $\nu$ be the unit normal on $\Gamma_H$ pointing into $U_H$.

Then, in the sense of tempered distributions in $t$, the following identity holds:
\begin{equation}\label{eq:tdhk}
\begin{aligned}
&
\int_{\Gamma_H}
\left(
\frac{\partial\Phi(x,t;y)}{\partial\nu(x)}
*
\Phi(x,-t;z)
-
\Phi(x,t;y)
*
\frac{\partial\Phi(x,-t;z)}{\partial\nu(x)}
\right)
\, \mathrm ds(x)
\\
&=
\frac1{4\pi^2}
\int_{S_+^1}
\operatorname{p.v.}
\frac1{
t-c^{-1}\hat x\cdot(z-y)
}
\, \mathrm ds(\hat x)
=K(t;y,z),
\end{aligned}
\end{equation}
where $\operatorname{p.v.}$ denotes the Cauchy principal value.
\end{lemma}

\begin{proof}
We first derive the corresponding frequency-domain identity and then apply the inverse Fourier transform.

For $f\in\mathcal S'(\mathbb R)$, define the Fourier transform by
\[
\widehat f(\omega)
=
\mathcal F[f](\omega)
=
\int_{-\infty}^{\infty}
f(t)e^{\mathrm i\omega t}\,dt.
\]

For $\omega\neq0$, let
\[
k=\frac{|\omega|}{c}.
\]
The Fourier transform of the causal wave fundamental solution satisfies
\[
\mathcal F[\Phi(x,\cdot;y)](\omega)
=
G_\omega(x,y),
\]
where $G_\omega(x,y)$ is the radiating fundamental solution of the Helmholtz equation
\[
(\Delta+\omega^2/c^2)G_\omega(x,y)
=
-\delta(x-y).
\]

For $\omega>0$,
\[
G_\omega(x,y)
=
\frac{\mathrm i}{4}H_0^{(1)}(k|x-y|),
\]
where $H_0^{(1)}$ is the Hankel function of the first kind of order zero.
Moreover,
\[
G_{-\omega}(x,z)
=
\overline{G_\omega(x,z)}.
\]

Let
\[
B_R^+
=
\{x\in U_H:|x|<R\}.
\]
Since $y,z\notin U_H$, both
\[
G_\omega(x,y)
\quad\mbox{and}\quad
G_{-\omega}(x,z)
\]
satisfy the homogeneous Helmholtz equation in $B_R^+$.
Applying Green's second identity in $B_R^+$ yields
\[
\int_{\partial B_R^+}
\left(
\frac{\partial G_\omega(x,y)}{\partial n}
G_{-\omega}(x,z)
-
G_\omega(x,y)
\frac{\partial G_{-\omega}(x,z)}{\partial n}
\right)
\, \mathrm ds(x)
=
0,
\]
where $n$ denotes the outward unit normal on $\partial B_R^+$.

The boundary $\partial B_R^+$ consists of the horizontal segment on $\Gamma_H$ and the upper semicircle $\partial B_R^+\cap U_H$.
Since the outward normal on $\Gamma_H$ equals $-\nu$, we obtain
\begin{equation}\label{eq:greenbd}
\begin{aligned}
&
\int_{\Gamma_H}
\left(
\frac{\partial G_\omega(x,y)}{\partial\nu(x)}
G_{-\omega}(x,z)
-
G_\omega(x,y)
\frac{\partial G_{-\omega}(x,z)}{\partial\nu(x)}
\right)
\, \mathrm ds(x)
\\
&=
\lim_{R\to\infty}
\int_{\partial B_R^+\cap U_H}
\left(
\frac{\partial G_\omega(x,y)}{\partial r}
G_{-\omega}(x,z)
-
G_\omega(x,y)
\frac{\partial G_{-\omega}(x,z)}{\partial r}
\right)
\, \mathrm ds(x).
\end{aligned}
\end{equation}

For $x=R\hat x$, $\hat x\in S_+^1$, the Helmholtz fundamental solution admits the asymptotic expansion
\[
G_\omega(R\hat x,y)
=
\frac{C_\omega}{\sqrt R}
\exp\left(
\mathrm i\,\operatorname{sgn}(\omega)kR
-
\mathrm i\frac{\omega}{c}\hat x\cdot y
\right)
+
O(R^{-3/2}),
\]
where
\[
|C_\omega|^2
=
\frac1{8\pi k}.
\]
Moreover,
\[
\frac{\partial G_\omega(R\hat x,y)}{\partial r}
=
\mathrm i\,\operatorname{sgn}(\omega)k
G_\omega(R\hat x,y)
+
O(R^{-3/2}),
\]
and similarly
\[
\frac{\partial G_{-\omega}(R\hat x,z)}{\partial r}
=
-
\mathrm i\,\operatorname{sgn}(\omega)k
G_{-\omega}(R\hat x,z)
+
O(R^{-3/2}).
\]

Substituting the above asymptotic expansions into \eqref{eq:greenbd} gives
\[
\frac{\partial G_\omega}{\partial r}G_{-\omega}
-
G_\omega\frac{\partial G_{-\omega}}{\partial r}
=
2 \mathrm i\,\operatorname{sgn}(\omega)k
G_\omega G_{-\omega}
+
O(R^{-2}).
\]

Furthermore,
\[
G_\omega(R\hat x,y)
G_{-\omega}(R\hat x,z)
=
\frac{|C_\omega|^2}{R}
\exp\left(
\mathrm i\frac{\omega}{c}\hat x\cdot(z-y)
\right)
+
O(R^{-2}).
\]

Since $\mathrm ds(x)=R\,\mathrm ds(\hat x)$ on the semicircle, we obtain
\[
\begin{aligned}
&
\lim_{R\to\infty}
\int_{\partial B_R^+\cap U_H}
\left(
\frac{\partial G_\omega}{\partial r}
G_{-\omega}
-
G_\omega
\frac{\partial G_{-\omega}}{\partial r}
\right)
\,\mathrm ds(x)
\\
&=
\frac{
\mathrm i\,\operatorname{sgn}(\omega)
}{4\pi}
\int_{S_+^1}
\exp\left(
\mathrm i\frac{\omega}{c}\hat x\cdot(z-y)
\right)
\,\mathrm ds(\hat x).
\end{aligned}
\]

Hence,
\begin{equation}\label{eq:fdhk}
\begin{aligned}
&
\int_{\Gamma_H}
\left(
\frac{\partial G_\omega(x,y)}{\partial\nu(x)}
G_{-\omega}(x,z)
-
G_\omega(x,y)
\frac{\partial G_{-\omega}(x,z)}{\partial\nu(x)}
\right)
\,\mathrm ds(x)
\\
&=
\frac{
\mathrm i\,\operatorname{sgn}(\omega)
}{4\pi}
\int_{S_+^1}
\exp\left(
\mathrm i\frac{\omega}{c}\hat x\cdot(z-y)
\right)
\,\mathrm ds(\hat x).
\end{aligned}
\end{equation}

We now return to the time domain.
Using
\[
\mathcal F_t[\Phi(x,t;y)]
=
G_\omega(x,y),
\]
and
\[
\mathcal F_t[\Phi(x,-t;z)]
=
G_{-\omega}(x,z),
\]
together with the convolution theorem, the inverse Fourier transform of the left-hand side of \eqref{eq:fdhk} equals
\[
\int_{\Gamma_H}
\left(
\frac{\partial\Phi(x,t;y)}{\partial\nu(x)}
*
\Phi(x,-t;z)
-
\Phi(x,t;y)
*
\frac{\partial\Phi(x,-t;z)}{\partial\nu(x)}
\right)
\, \mathrm ds(x).
\]

Finally, using the distributional identity
\[
\mathcal F_t^{-1}
\left[
\mathrm i\,\operatorname{sgn}(\omega)e^{\mathrm i\omega\tau}
\right]
=
\operatorname{p.v.}
\frac1{\pi(t-\tau)},
\]
with
\[
\tau
=
c^{-1}\hat x\cdot(z-y),
\]
we conclude that
\[
\mathcal F_t^{-1}
\left[
\frac{
\mathrm i\,\operatorname{sgn}(\omega)
}{4\pi}
e^{\mathrm i\omega\tau}
\right]
=
\frac1{4\pi^2}
\operatorname{p.v.}
\frac1{t-\tau}.
\]

This proves \eqref{eq:tdhk}.
\end{proof}

\begin{lemma}[Pulse-filtered Helmholtz--Kirchhoff identity]
For $y,z\in \mathbb R^2\setminus \overline{U_H}$, the pulse-filtered fundamental solutions satisfy
\[
\begin{aligned}
&
\int_{\Gamma_H}
\left(
\partial_\nu\Phi_\chi(x,t;y)
*
\Phi_\chi(x,-t;z)
-
\Phi_\chi(x,t;y)
*
\partial_\nu\Phi_\chi(x,-t;z)
\right)
\,\mathrm ds(x)
\\
&=
\big(\chi*K(\cdot;y,z)*\chi(-\cdot)\big)(t).
\end{aligned}
\]
\end{lemma}

The preceding identities are understood in the sense of tempered distributions in time. The kernel $K(t;y,z)$ is a half-aperture kernel because only the upper measurement line $\Gamma_H$ is used.

\subsection{Back-propagation and virtual scattering field}

We now use the measured pulse-filtered Cauchy data
\[
\{v_\chi^s(x,t;y),\partial_\nu v_\chi^s(x,t;y):x\in\Gamma_H\}
\]
to define a time-domain back-propagation field. For a sampling point $z\in D^+$, set
\[
V^i(y,z,t)
:=
\big(\chi*K(\cdot;y,z)*\chi(-\cdot)\big)(t),
\]
which will play the role of a virtual incident field.

Motivated by Lemma~4.1 and the pulse-filtered identity above, we define the back-propagated field associated with the sampling point $z\in D^+$ by
\begin{equation}\label{eq:backprop}
\begin{aligned}
V^s(y,z,t)
:=
\int_{\Gamma_H}
\Big(
&
\partial_\nu v_\chi^s(x,t;y)
*
\Phi_\chi(x,-t;z)
\\
&
-
v_\chi^s(x,t;y)
*
\partial_\nu\Phi_\chi(x,-t;z)
\Big)
\,\mathrm ds(x).
\end{aligned}
\end{equation}

The above definition propagates the measured Cauchy data on $\Gamma_H$ back into the sampling region. The following theorem is the main result of this section: it identifies the back-propagated field as a genuine scattered field generated by the virtual incident field $V^i$.

\begin{theorem}\label{thm:pre_imagingFunction}
Let $z\in\mathbb{R}^2\backslash\overline{U_H}$ be a sampling point, and let $V^s(y,z,t)$ be defined by \eqref{eq:backprop}. Then $V^s(y,z,t)$ is the scattered field generated by the virtual incident field
\[
V^i(y,z,t)
=
\big(\chi*K(\cdot;y,z)*\chi(-\cdot)\big)(t).
\]
In particular,
\[
c^{-2}\partial_{tt}V^s(y,z,t)
-
\Delta_yV^s(y,z,t)
=
0,
\qquad y\in D^+,
\]
and on the rough surface,
\[
V^s(y,z,t)
=
-
V^i(y,z,t),
\qquad y\in\Gamma.
\]
\end{theorem}
\begin{proof}
For fixed $x\in\Gamma_H$ and $z\in D^+$, the functions
\[
\Phi_\chi(x,-t;z),
\qquad
\partial_\nu\Phi_\chi(x,-t;z)
\]
do not depend on the variable $y$.
Using the reciprocity relation
\[
v_\chi^s(x,t;y)
=
v_\chi^s(y,t;x),
\]
we may regard the field $v_\chi^s(x,t;y)$ as a solution of the homogeneous wave equation with respect to the variable $y$.
Consequently, the back-propagated field $V^s(y,z,t)$ defined by \eqref{eq:backprop} satisfies
\[
c^{-2}\partial_{tt}V^s(y,z,t)
-
\Delta_yV^s(y,z,t)
=
0,
\qquad y\in D^+.
\]

Next, we identify the boundary trace of $V^s$ on the rough surface $\Gamma$.
For $y\in\Gamma$, the Dirichlet boundary condition for the pulse-filtered scattered field yields
\[
v_\chi^s(y,t;x)
=
-
\Phi_\chi(y,t;x).
\]
Applying the reciprocity relation again gives
\[
v_\chi^s(x,t;y)
=
-
\Phi_\chi(x,t;y),
\qquad y\in\Gamma.
\]

Substituting the above identity into the definition of $V^s$, we obtain
\[
\begin{aligned}
V^s(y,z,t)
=
-
\int_{\Gamma_H}
\Big(
&
\partial_\nu\Phi_\chi(x,t;y)
*
\Phi_\chi(x,-t;z)
\\
&
-
\Phi_\chi(x,t;y)
*
\partial_\nu\Phi_\chi(x,-t;z)
\Big)
\,\mathrm ds(x),
\qquad y\in\Gamma.
\end{aligned}
\]

By the pulse-filtered Helmholtz--Kirchhoff identity, the integral on the right-hand side equals
\[
(\chi*K(\cdot;y,z)*\chi(-\cdot))(t).
\]
Therefore,
\[
V^s(y,z,t)
=
-
V^i(y,z,t),
\qquad y\in\Gamma.
\]

Hence, $V^s$ satisfies the homogeneous wave equation in $D^+$ together with the Dirichlet boundary condition
\[
V^s=-V^i
\qquad \mbox{on } \Gamma.
\]
Since $V^s$ is causal in time, it satisfies the same initial-boundary value problem as the scattered field generated by the virtual incident field $V^i$.
By uniqueness of the time-domain scattering problem, we conclude that $V^s$ coincides with the scattered field corresponding to the virtual incident wave $V^i$.
\end{proof}

\begin{remark}
[Well-definedness of the back-propagation operator] The back-propagated field $V^s(y,z,t)$ is well-defined in the sense of tempered distributions with respect to the time variable. Indeed, since the pulse function $\chi$ is causal and sufficiently smooth with compact or rapidly decaying temporal support, both $v_\chi^s$ and $\Phi_\chi$ are tempered distributions in time and are smooth away from their source points. Moreover, for fixed $y,z\in D^+$ and finite observation time $t$, the finite propagation speed of the wave equation implies that only a bounded portion of the measurement line $\Gamma_H$ contributes to the convolution integral in \eqref{eq:backprop}. More precisely, when $|x_1|$ is sufficiently large, the travel time from the source point to the measurement point $x$, and then from $x$ back to the sampling point $z$, exceeds the observation time window under consideration. Consequently, the corresponding contribution vanishes. Therefore, although the integration is taken over the infinite line $\Gamma_H$, the effective support of the integrand is bounded for finite observation times. If the pulse function is not compactly supported but rapidly decaying, the same conclusion follows from the far-field decay of the two-dimensional wave field together with the rapid temporal decay of $\chi$.
\end{remark}

\begin{remark}
[Reciprocity relation] The reciprocity relation \[ v_\chi^s(x,t;y) = v_\chi^s(y,t;x) \] used in the proof of Theorem~\ref{thm:pre_imagingFunction} follows from the symmetry of the time-domain Dirichlet scattering problem in a homogeneous acoustic medium. More precisely, the wave operator is self-adjoint, while the Dirichlet boundary condition on the rough surface is symmetric with respect to the source and observation variables. Let $v_\chi^s(\cdot,t;y)$ and $v_\chi^s(\cdot,t;x)$ denote the scattered fields generated by pulse-filtered point sources located at $y$ and $x$, respectively. Applying Green's second identity in a truncated domain, using the Dirichlet boundary condition on the rough surface, and then passing to the limit with the radiation condition yield the source--receiver reciprocity relation above. The same argument also applies to the corresponding normal derivatives on the measurement line $\Gamma_H$, which justifies the interchange of source and observation variables in the back-propagation formula.
\end{remark}

\subsection{Focusing property of the full-aperture kernel}

The kernel obtained from Lemma~4.1 involves only the upper half aperture $S_+^1$, since the Cauchy data are measured on the horizontal line $\Gamma_H$. To compensate for the missing lower-aperture contribution, we introduce the full-aperture kernel and a corresponding correction kernel. Let
\[
S_-^1
:=
\{\hat x\in\mathbb R^2:|\hat x|=1,\ \hat x_2<0\},
\qquad
S^1=S_+^1\cup S_-^1,
\]
and define the full-aperture Helmholtz--Kirchhoff kernel by
\begin{equation}\label{eq:Kfull}
K_{\rm full}(t;y,z)
:=
\frac1{4\pi^2}
\int_{S^1}
\operatorname{p.v.}
\frac1{
t-c^{-1}\hat x\cdot(z-y)
}
\,\mathrm ds(\hat x).
\end{equation}
Here $y$ denotes the source or observation variable and $z\in D^+$ is a sampling point. If $K(t;y,z)$ denotes the half-aperture kernel in Lemma~4.1, we set
\begin{equation}\label{eq:Kcorr}
K_{\rm corr}(t;y,z)
:=
K_{\rm full}(t;y,z)-K(t;y,z).
\end{equation}
The correction kernel $K_{\rm corr}$ represents a synthetic lower-aperture contribution. It does not introduce additional measured data; it compensates for the aperture missing from the measured-line Helmholtz--Kirchhoff identity.

\begin{remark}
[Qualitative focusing property of the full-aperture kernel]
The full-aperture kernel $K_{\rm full}(t;y,z)$ is rotationally invariant in the spatial difference $z-y$ and therefore depends on $y$ and $z$ only through $r=|z-y|$. As $r$ decreases, the main features of the regularized full-aperture profile become increasingly concentrated near $t=0$. In this qualitative sense, $K_{\rm full}$ provides a focusing kernel for sampling points close to the boundary.
Indeed, the integral in \eqref{eq:Kfull} is taken over the full unit circle $S^1$, so it is invariant under rotations of the spatial difference $z-y$. Moreover, the delay term $c^{-1}\hat x\cdot(z-y)$ ranges over values whose magnitude is bounded by $r/c$. Thus, when $r$ is small, the relevant temporal contributions of the regularized or pulse-filtered kernel are concentrated near $t=0$.
\end{remark}

\begin{figure}[htbp]
    \centering
    \begin{subfigure}[t]{0.48\linewidth}
        \centering
        \includegraphics[width=\linewidth]{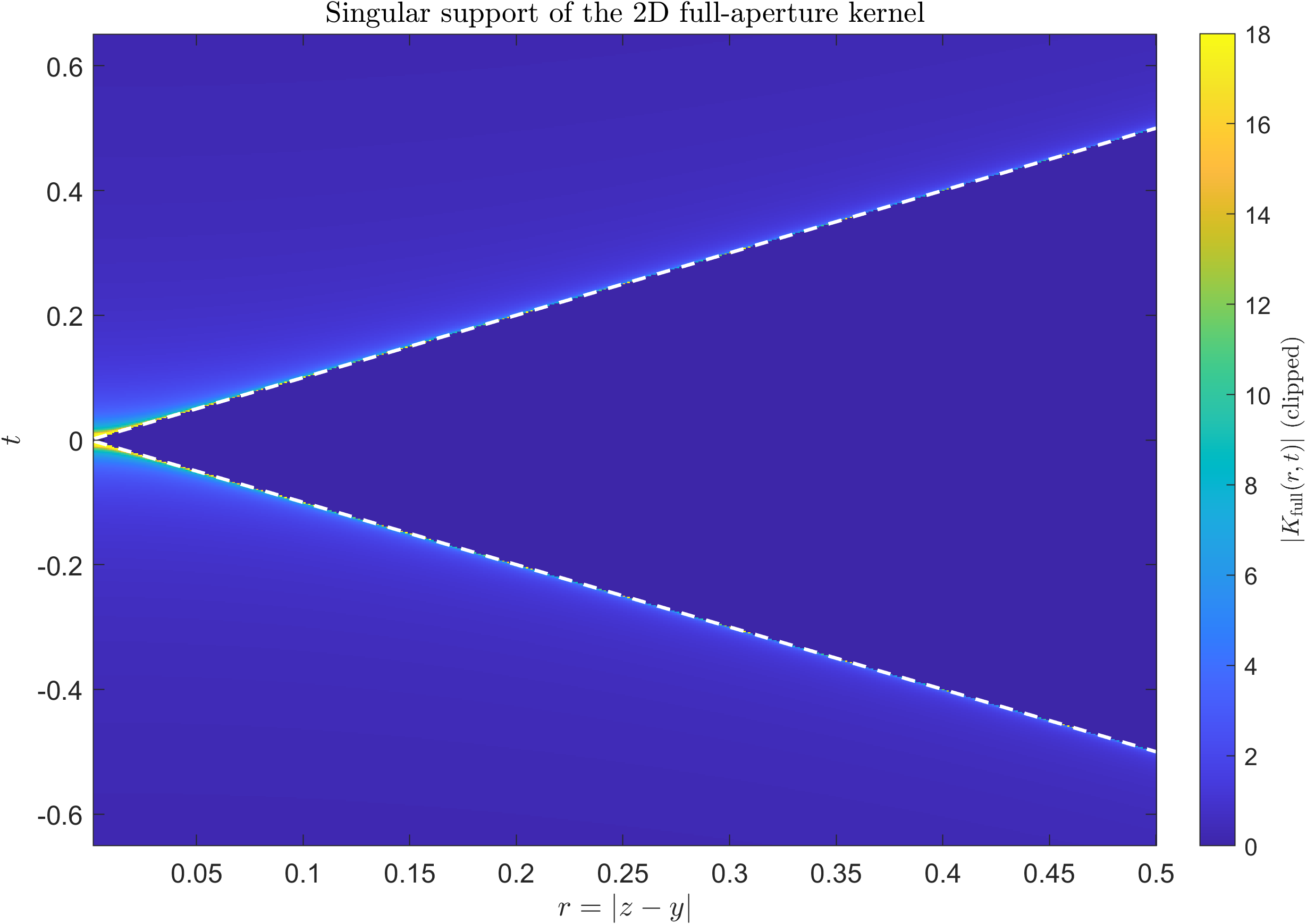}
        \caption{Regularized magnitude of $K_{\rm full}(r,t)$.}
        \label{fig:kfull_heatmap}
    \end{subfigure}
    \hfill
    \begin{subfigure}[t]{0.48\linewidth}
        \centering
        \includegraphics[width=\linewidth]{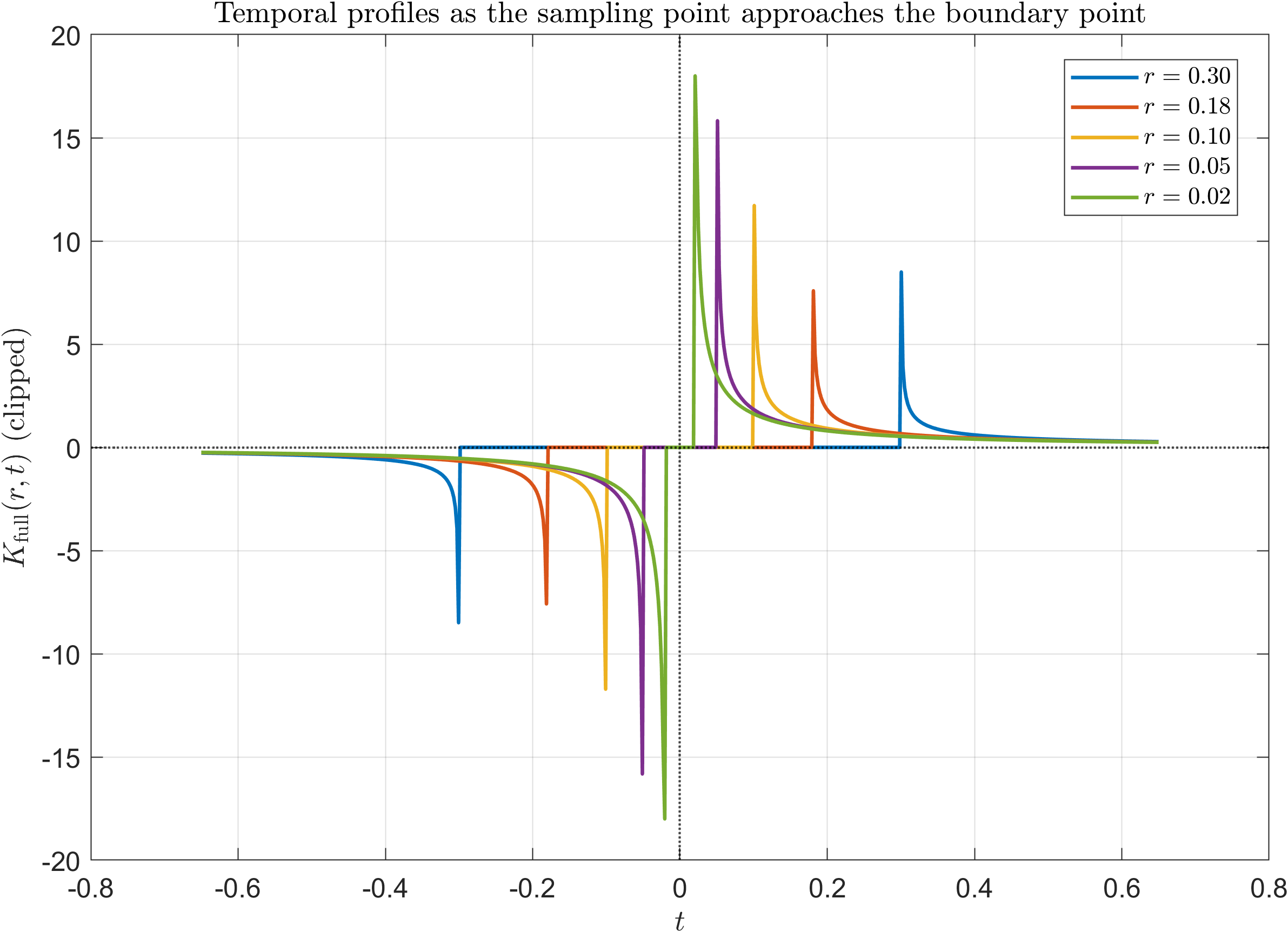}
        \caption{Temporal profiles for decreasing values of $r$.}
        \label{fig:kfull_profiles}
    \end{subfigure}
    \caption{Singular focusing behavior of the two-dimensional full-aperture Helmholtz--Kirchhoff kernel.  The dashed curves in (a) indicate the characteristic singularities $t=\pm r/c$.  As the sampling point $z$ approaches the boundary point $y$, namely $r=|z-y|\to0$, these singularities move toward $t=0$ and the regularized kernel becomes increasingly concentrated.}
    \label{fig:boundarydata_behavior}
\end{figure}

Figure~\ref{fig:boundarydata_behavior} illustrates this focusing structure in a regularized numerical sense.

\begin{remark}
There is an important difference between the two- and three-dimensional time-domain kernels. In two dimensions, the Helmholtz--Kirchhoff kernel is a principal-value distribution and the wave equation does not satisfy the sharp Huygens principle. Hence the focusing profile should be interpreted as a regularized or pulse-filtered concentration effect rather than a delta-type peak. Nevertheless, the full-aperture kernel \eqref{eq:Kfull}, together with the pulse filtering used in the back-propagation field, still retains sufficient local concentration to support the direct sampling method.
\end{remark}

\subsection{Imaging functional and imaging mechanism}

The correction field associated with the missing lower-aperture contribution is defined by
\begin{equation}\label{eq:Vcorr}
V_{\rm corr}^i(y,z,t)
:=
\big(\chi*K_{\rm corr}(\cdot;y,z)*\chi(-\cdot)\big)(t).
\end{equation}
The field $V_{\rm corr}^i$ should be understood as a synthetic contribution associated with the lower half aperture $S_-^1$. Its role is not to introduce additional measured data, but to compensate for the focusing structure suggested by the full-aperture Helmholtz--Kirchhoff identity.

For a finite observation time $T>0$, we define the imaging functional by
\begin{equation}\label{eq:imagingFunction}
I(z)
:=
\int_{\Gamma_H}\int_0^T
\left|
V^s(y,z,t)+V_{\rm corr}^i(y,z,t)
\right|^2
\,\mathrm{d}t\,\mathrm{d}s(y),
\qquad z\in D^+.
\end{equation}
Thus $I(z)$ measures the time-integrated energy of the corrected back-propagated field associated with the sampling point $z$.

The virtual scattering theorem shows that $V^s$ behaves as the scattered field generated by the virtual incident field associated with the half-aperture kernel $K$. The correction field $V_{\rm corr}^i$ adds the synthetic lower-aperture contribution associated with $K_{\rm corr}$. Consequently, the corrected field $V^s+V_{\rm corr}^i$ is designed to reflect the focusing behavior of the full-aperture kernel $K_{\rm full}$. When the sampling point $z$ approaches the rough surface $\Gamma$, the virtual incident field interacts more strongly with the Dirichlet boundary, producing larger corrected back-propagated energy. When $z$ is away from $\Gamma$, contributions from different directions tend to cancel, and the energy is expected to be smaller. Thus local peaks of $I(z)$ provide an indicator for the rough surface.

\begin{remark}
The preceding mechanism is qualitative and is intended to explain the sampling behavior of \eqref{eq:imagingFunction}. It should not be interpreted as a rigorous resolution theorem without additional quantitative estimates.
\end{remark}

\providecommand{\resultplaceholder}[1]{%
    \fbox{\parbox[c][0.62\textwidth][c]{0.92\textwidth}{\centering\small #1}}}
\providecommand{\resultpanel}[2]{%
    \begin{minipage}{0.30\textwidth}
        \centering
        \includegraphics[width=\textwidth]{#1}
        {\resultplaceholder{Result to be inserted}}
        \subcaption{#2}
    \end{minipage}}
\providecommand{\resultpanelwide}[2]{%
    \begin{minipage}{0.45\textwidth}
        \centering
        \includegraphics[width=\textwidth]{#1}
        {\resultplaceholder{Result to be inserted}}
        \subcaption{#2}
    \end{minipage}}
    
\section{Numerical experiments}
\label{sec:numerical-experiments}

This section evaluates the proposed time-domain direct sampling method (TD-DSM).
We first validate the method on six representative rough surfaces and then compare it with single-frequency DSM and multi-frequency DSM.
We subsequently investigate resolution, measurement configuration, sampling density, and robustness with respect to noise.

Unless otherwise stated, the computational domain is $[-10,10]\times[-5,5]$, the observation time interval is $[0,5]$, and the sampling region is $\Omega=[-3,3]\times[-0.4,0.4]$.
Sources and receivers are uniformly distributed on the
truncated measurement line
$\Gamma_H^A=\{(x_1,H):x_1\in[-A,A]\}$.
Here, $N_s$ and $N_r$ denote the numbers of sources and
receivers, respectively, and $N_t$ denotes the number
of temporal sampling points.
The baseline parameters are $A=4$, $H=0.8$,
$N_s=N_r=30$, and $N_t=50$.
We use the causal pulse 
\begin{equation*}
\chi(t;k_c)=
\begin{cases}
\sin(k_ct)e^{-0.64(t-2)^2}, & t>0,\\
0, & t\leq 0,
\end{cases}
\end{equation*}
with center frequency $k_c=10$.
The baseline data are noise-free.

\begin{figure}[htbp]
    \centering
    \resultpanel{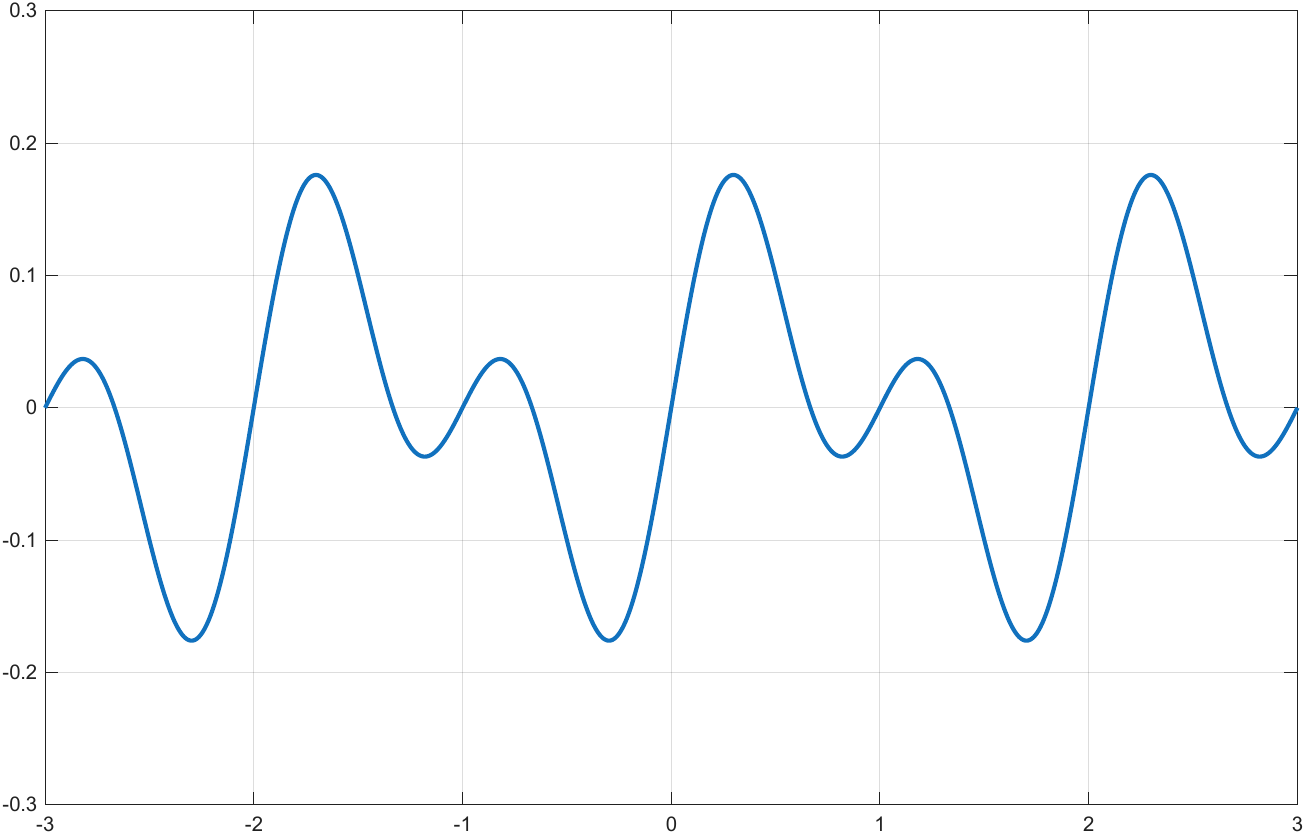}{$\Gamma_1$}\hfill
    \resultpanel{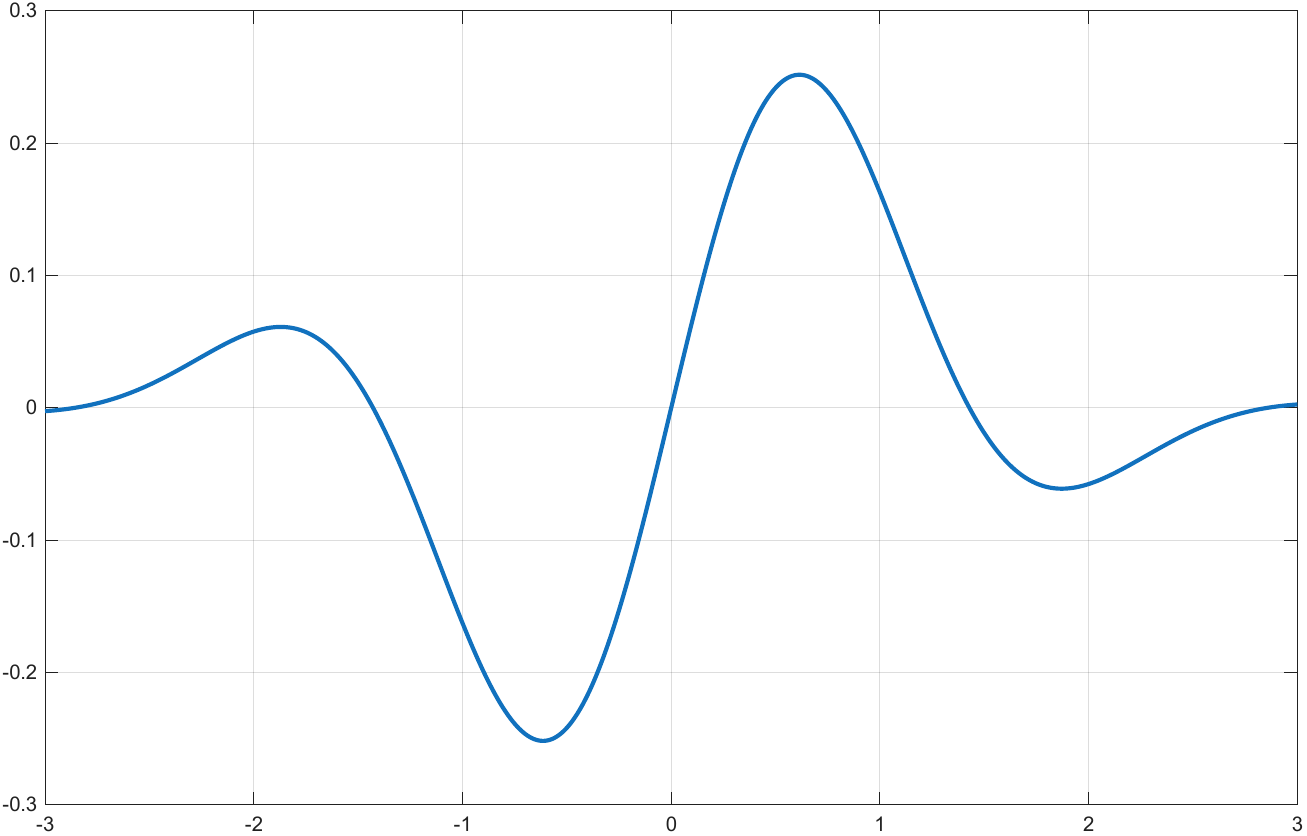}{$\Gamma_2$}\hfill
    \resultpanel{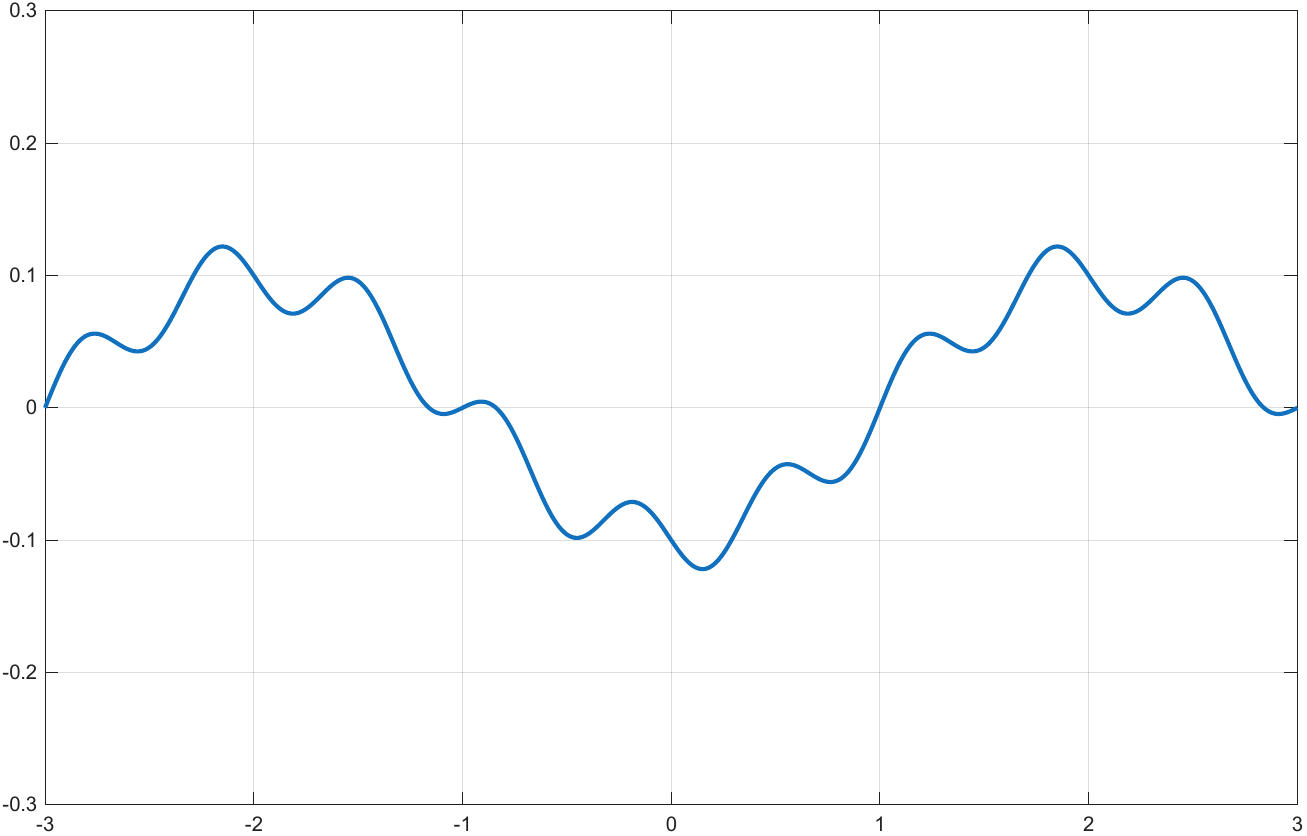}{$\Gamma_3$}

    \medskip
    \resultpanel{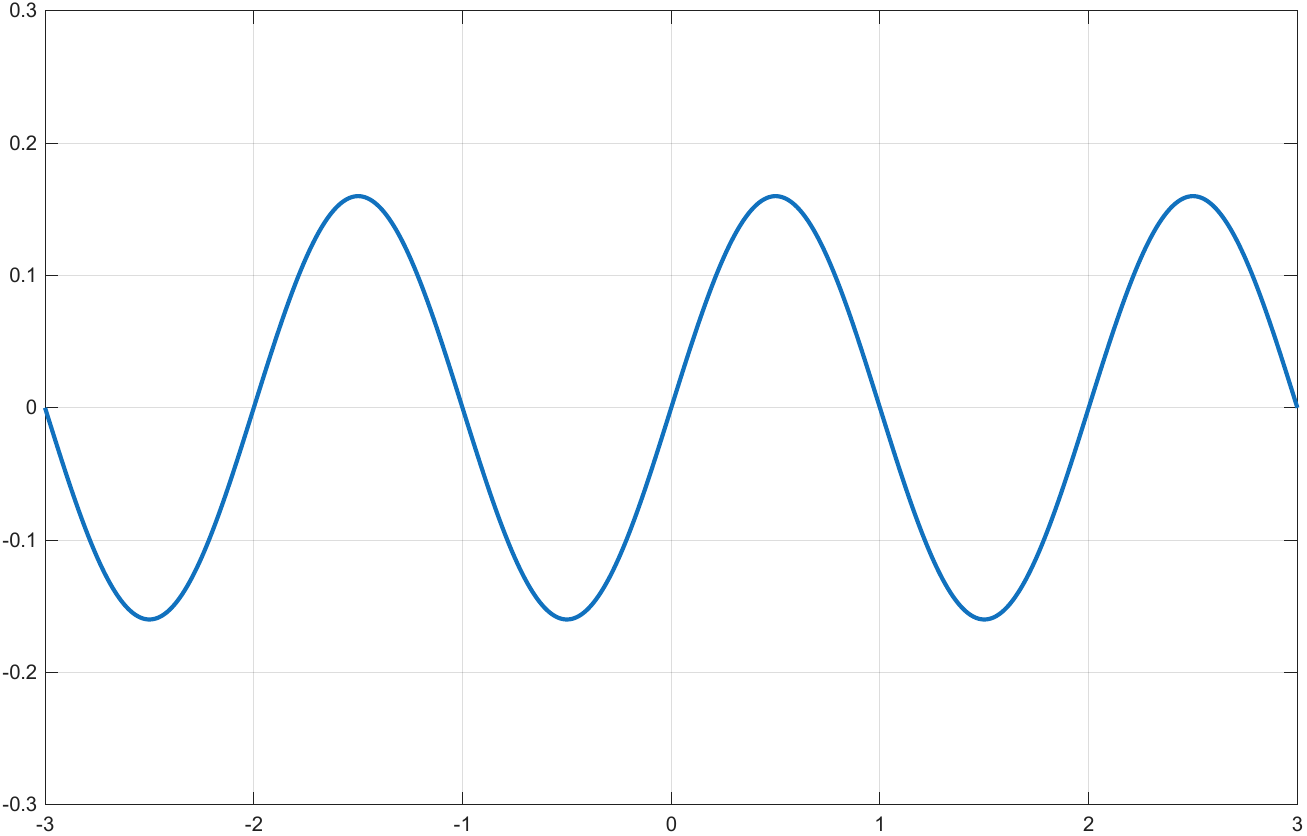}{$\Gamma_4$}\hfill
    \resultpanel{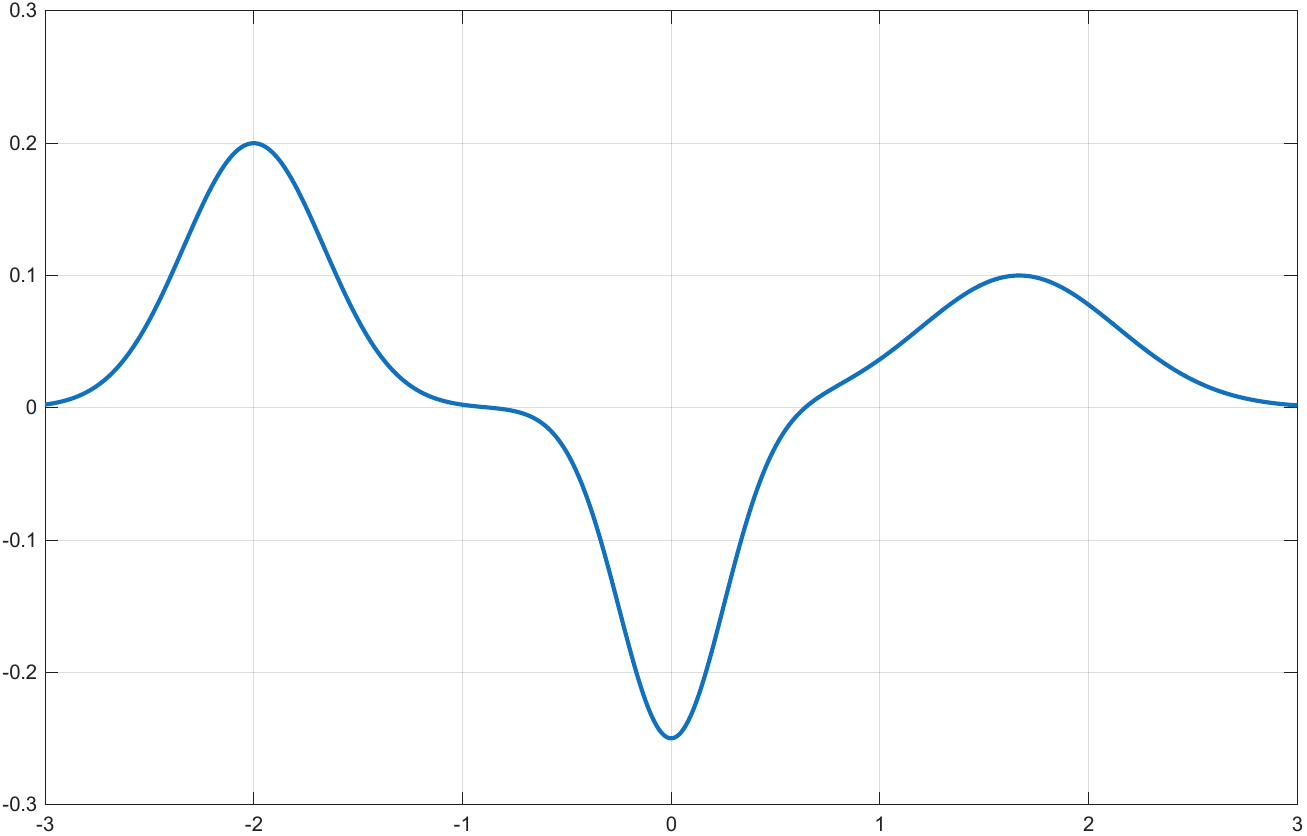}{$\Gamma_5$}\hfill
    \resultpanel{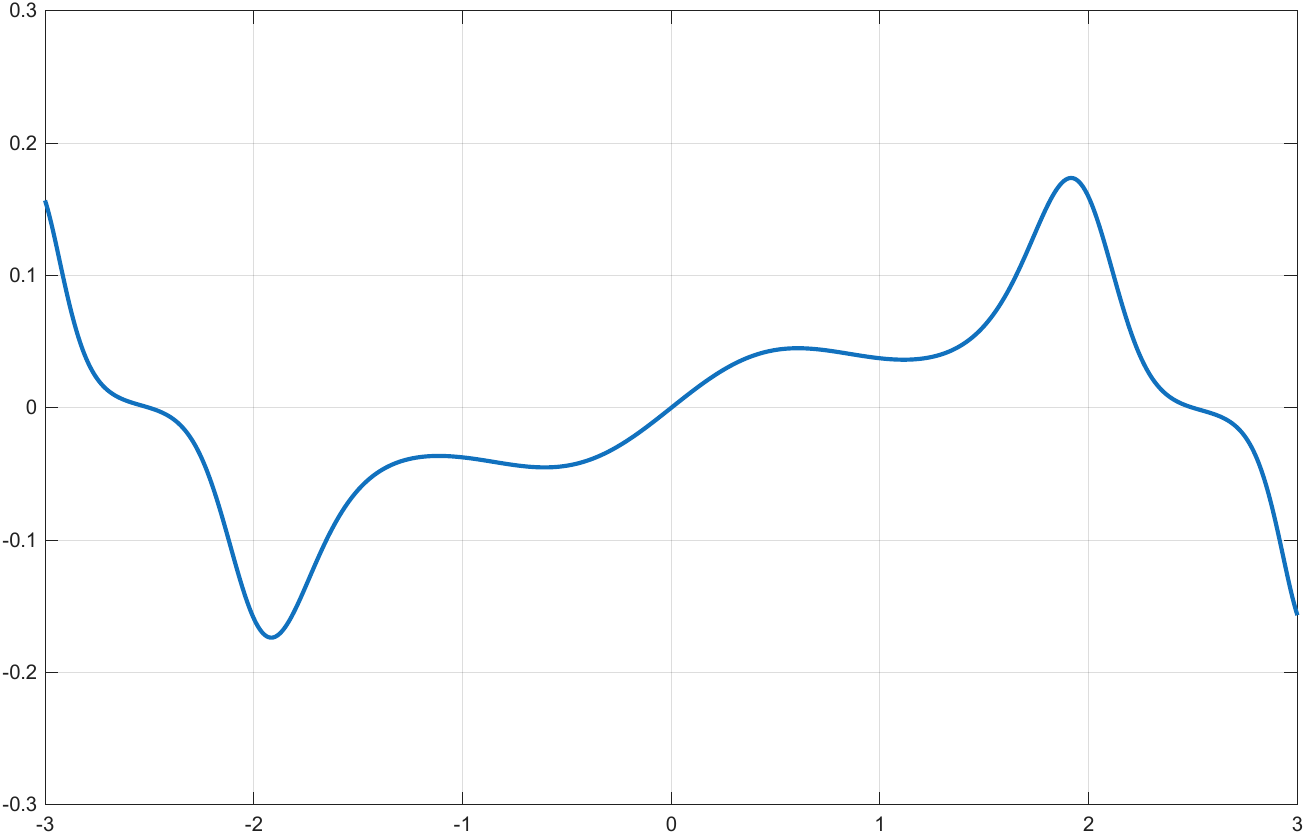}{$\Gamma_6$}
    \caption{Six representative rough surfaces used throughout the numerical
    experiments.}
    \label{fig:test-geometries}
\end{figure}

The six test surfaces are
\begin{align*}
    \Gamma_1&=\{(x_1,x_2):x_2=0.1\sin(2\pi x_1)+0.1\sin(\pi x_1)\},\\
    \Gamma_2&=\{(x_1,x_2):x_2=0.3\sin(0.7\pi x_1)e^{-0.4x_1^2}\},\\
    \Gamma_3&=\{(x_1,x_2):x_2=0.025\sin(3\pi(x_1-1))+0.1\sin(0.5\pi(x_1-1))\},\\
    \Gamma_4&=\{(x_1,x_2):x_2=0.16\sin(\pi x_1)\},\\
    \Gamma_5&=\{(x_1,x_2):x_2=0.1e^{-25(0.3x_1-0.5)^2}+0.2e^{-49(0.3x_1+0.6)^2}-0.25e^{-8x_1^2}\},\\
    \Gamma_6&=\{(x_1,x_2):x_2=0.1\sin(0.4\pi x_1)e^{-\sin(1.2x_1^2)}\},\\
\end{align*}
where $x_1\in\mathbb{R}$ in every definition, and their geometry is given in Figure~\ref{fig:test-geometries}.

\setcounter{roughpanel}{0}

For the numerical implementation, let $y_j$, $j=1,\ldots,N_s$, and $x_r$, $r=1,\ldots,N_r$, denote the uniformly spaced source and receiver locations on $\Gamma_H^A$, respectively, and let $t_n=(n-1)\Delta t$, $n=1,\ldots,N_t$, with $\Delta t=5/(N_t-1)$. The measurement line is truncated to $\Gamma_H^A$, and the spatial and temporal integrals in the back-propagation step are approximated by uniform quadrature rules.
In particular, the time convolution is approximated by
\[
    \big(a(\cdot)*b(\cdot)\big)(t_n)\approx \Delta t\sum_{q=1}^{n}a(t_q)b(t_{n-q+1}),
\]
and
\[
    \big(a(\cdot)*b(-\cdot)\big)(t_n)\approx \Delta t\sum_{q=n}^{N_t}a(t_q)b(t_{q-n+1}).
\]

For every source $y_j$ and sampling point
$z\in\Omega$, we approximate the corrected field
$V^s+V_{\mathrm{corr}}^i$ defined in Section~4
by the following discrete expression:
\begin{equation}\label{eq:discrete_backpropagation}
\begin{aligned}
\widetilde V(y_j,z,t_i)
={}&\Delta x\Delta t
\sum_{r=1}^{N_r}\sum_{l=i}^{N_t}
\Big[
\partial_\nu v_\chi^s(x_r,t_l;y_j)
\Phi_\chi(x_r,t_{l-i+1};z)\\
&\hspace{28mm}
-v_\chi^s(x_r,t_l;y_j)
\partial_\nu\Phi_\chi(x_r,t_{l-i+1};z)
\Big]\\
&-\Delta x\Delta t
\sum_{r=1}^{N_r}\sum_{l=i}^{N_t}
\Big[
\partial_\nu\Phi_\chi(x_r',t_l;y_j)
\Phi_\chi(x_r',t_{l-i+1};z)\\
&\hspace{28mm}
-\Phi_\chi(x_r',t_l;y_j)
\partial_\nu\Phi_\chi(x_r',t_{l-i+1};z)
\Big].
\end{aligned}
\end{equation}
Here $x_r'=(x_{r,1},-H)$ is the reflection of the receiver $x_r=(x_{r,1},H)$ across $x_2=0$. The normal derivatives in both terms are taken with respect to the upward normal $\nu=(0,1)$. The second term represents the lower-aperture correction through the reflected measurement line. The discrete imaging functional and its normalized version are
\begin{equation}
    \widetilde I(z)=\Delta y\Delta t
    \sum_{j=1}^{N_s}\sum_{i=1}^{N_t}
    \left|\widetilde V(y_j,z,t_i)\right|^2,
    \qquad
    I_{\mathrm{norm}}(z)=\frac{\widetilde I(z)}{\max_{z'\in\Omega}\widetilde I(z')}.
    \label{eq:discrete-indicator}
\end{equation}
The normalized indicator $I_{\mathrm{norm}}$ is used in all imaging maps below.

For quantitative comparison, we extract a ridge from the normalized indicator \(I_{\mathrm{norm}}\) and compute two metrics against the binary surface mask \(f_{\rm true}\).
A bidirectional local-window search first traces a continuous peak row index \(r_j\) along the columns of \(I_{\mathrm{norm}}\), and a Gaussian weight
\[
    W_{ij}=\exp\!\left(-\frac{(i-r_j)^2}{\sigma_w^2}\right)
\]
is then applied to give the weighted indicator \(I_W=I_{\mathrm{norm}}\odot W\), where $\odot$ is the element-wise product operator.
We report
\[
    E_{\rm rel}=
    \frac{\|I_W-f_{\rm true}\|_{L^2}}{\|I_W\|_{L^2}},
    \qquad
    C=
    \frac{\operatorname{mean}\{I_{W,ij}:f_{\rm true,ij}=1\}}
         {\operatorname{mean}\{I_{W,ij}:f_{\rm true,ij}=0\}},
\]
where \(\|\cdot\|_{L^2}\) is understood discretely on the sampling grid. The same ridge-extraction rule and weighting parameters are used in all experiments.

\subsection{Validation of the proposed method}
\label{subsec:validation}

We first apply the proposed TD-DSM to the six rough surfaces using the default configuration. These examples include periodic oscillations, localized perturbations, and more complicated surface profiles, allowing us to examine the applicability of the method to different geometric features.

\begin{figure}[htbp]
    \centering
    \resultpanel{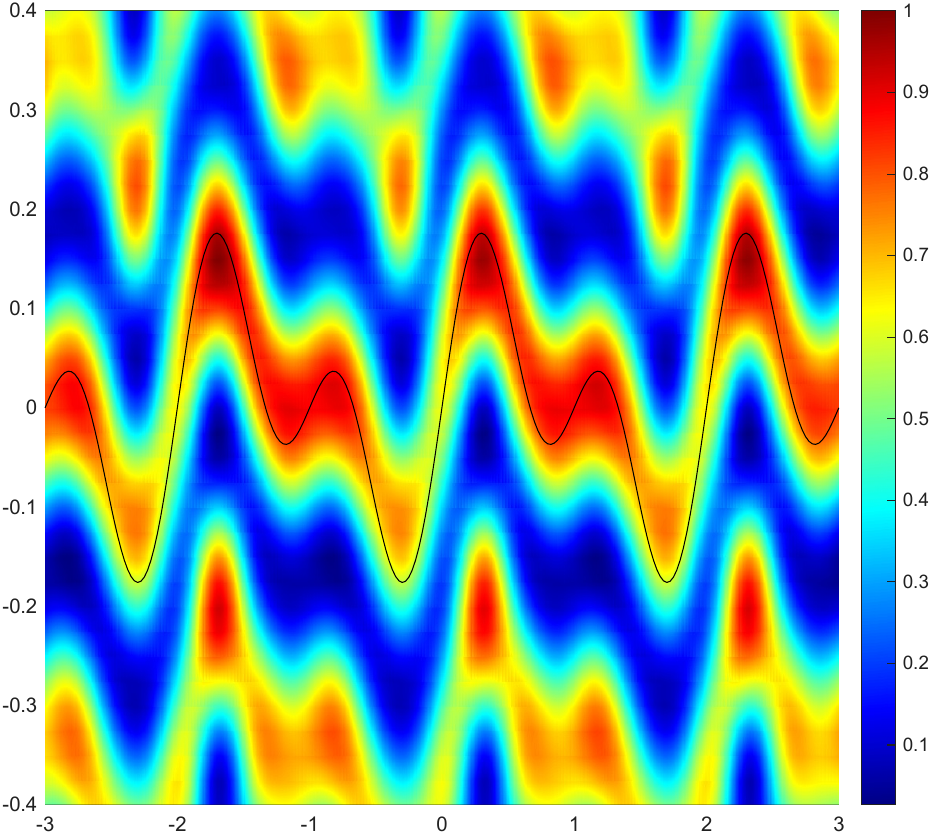}{$\Gamma_1$}\hfill
    \resultpanel{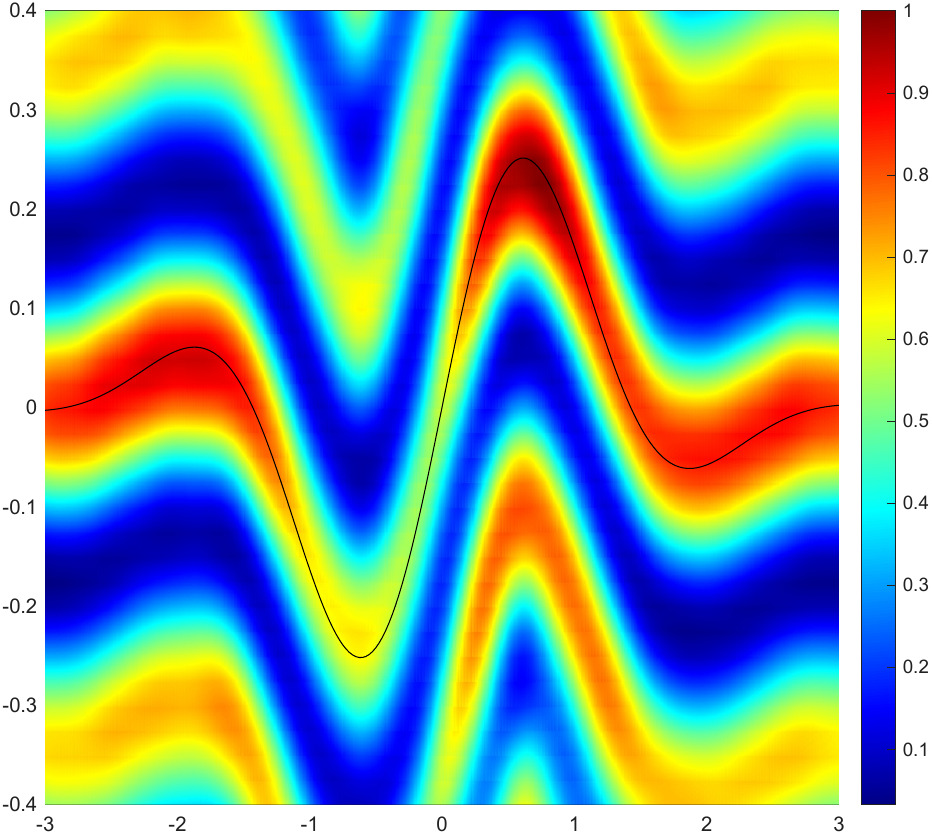}{$\Gamma_2$}\hfill
    \resultpanel{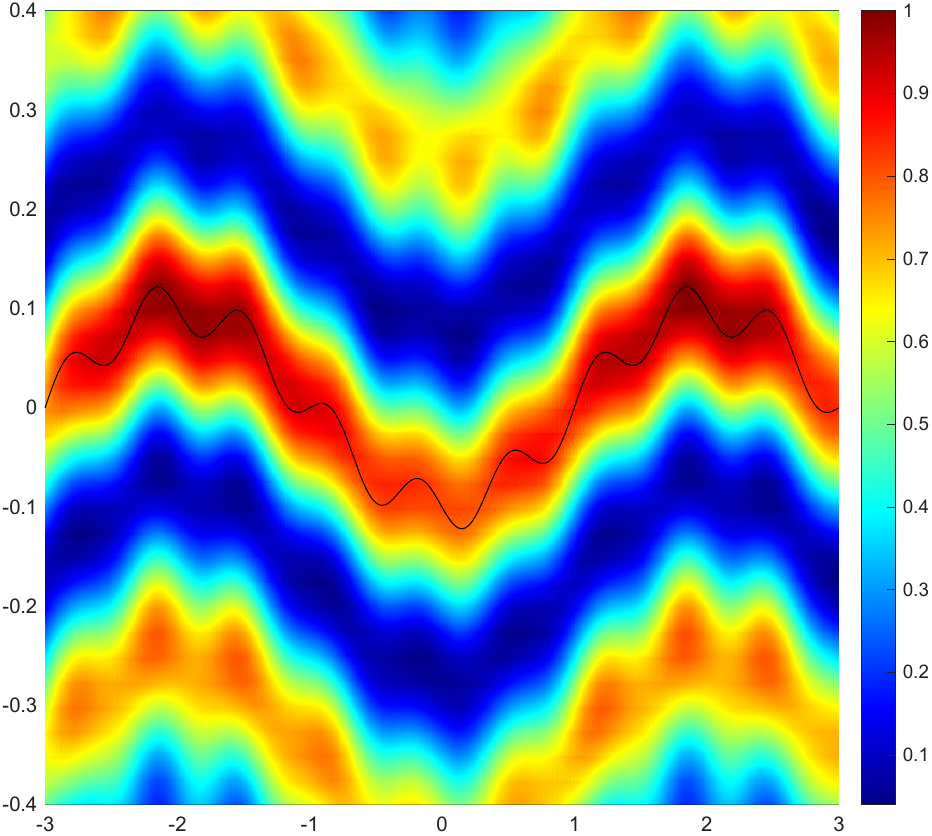}{$\Gamma_3$}

    \medskip
    \resultpanel{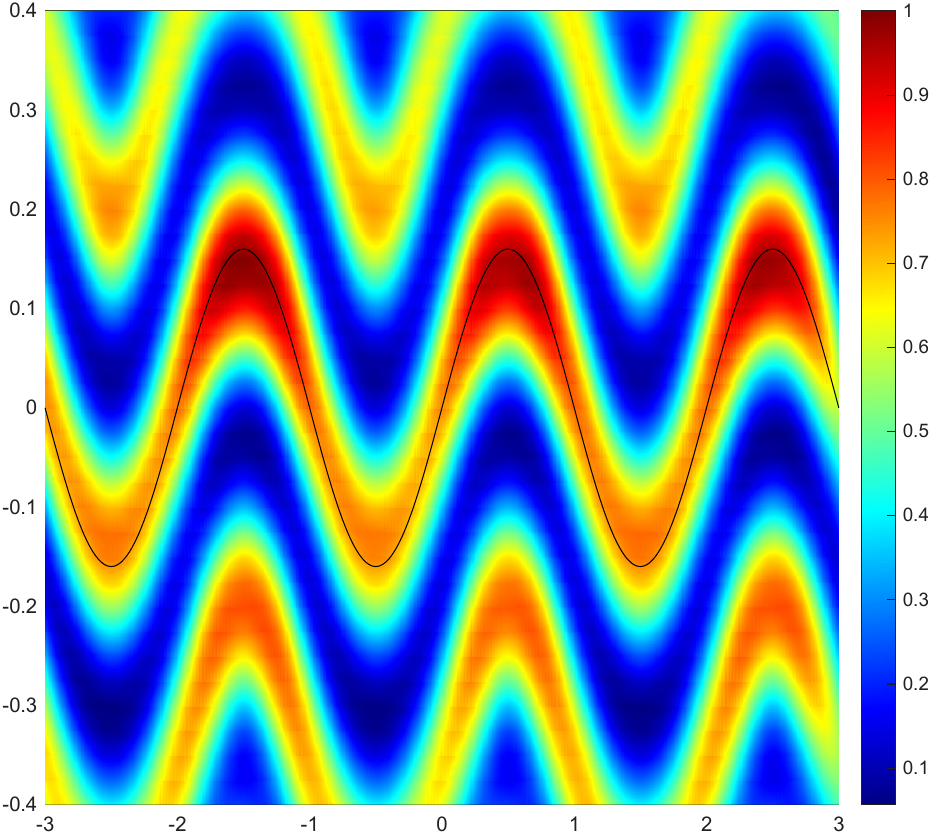}{$\Gamma_4$}\hfill
    \resultpanel{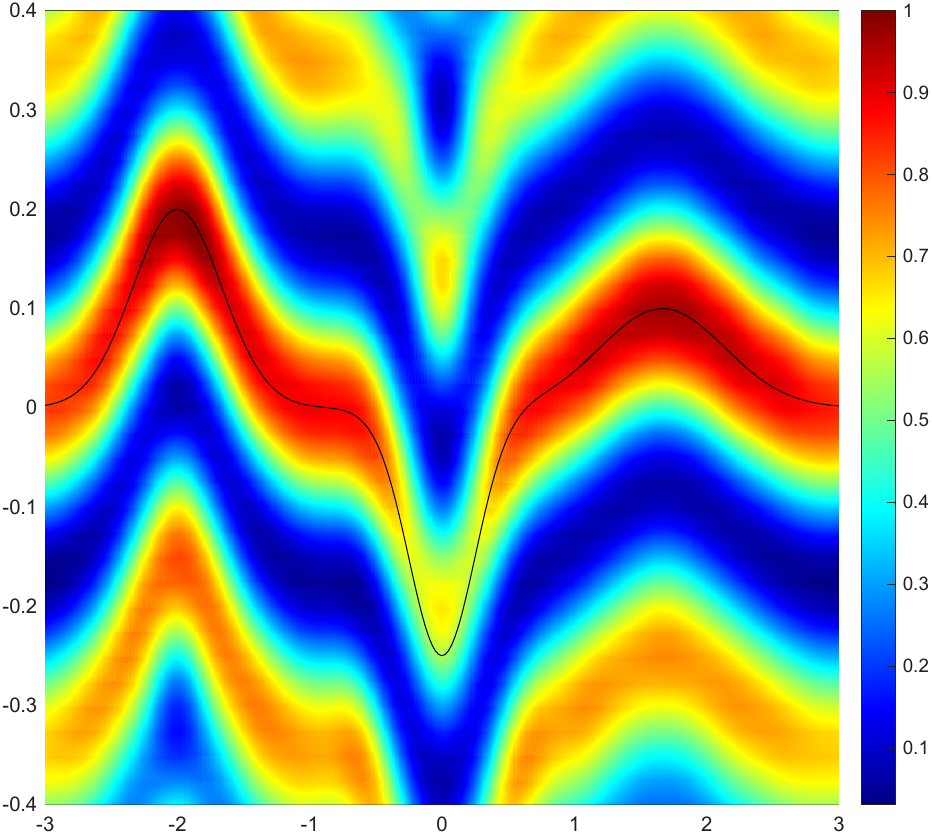}{$\Gamma_5$}\hfill
    \resultpanel{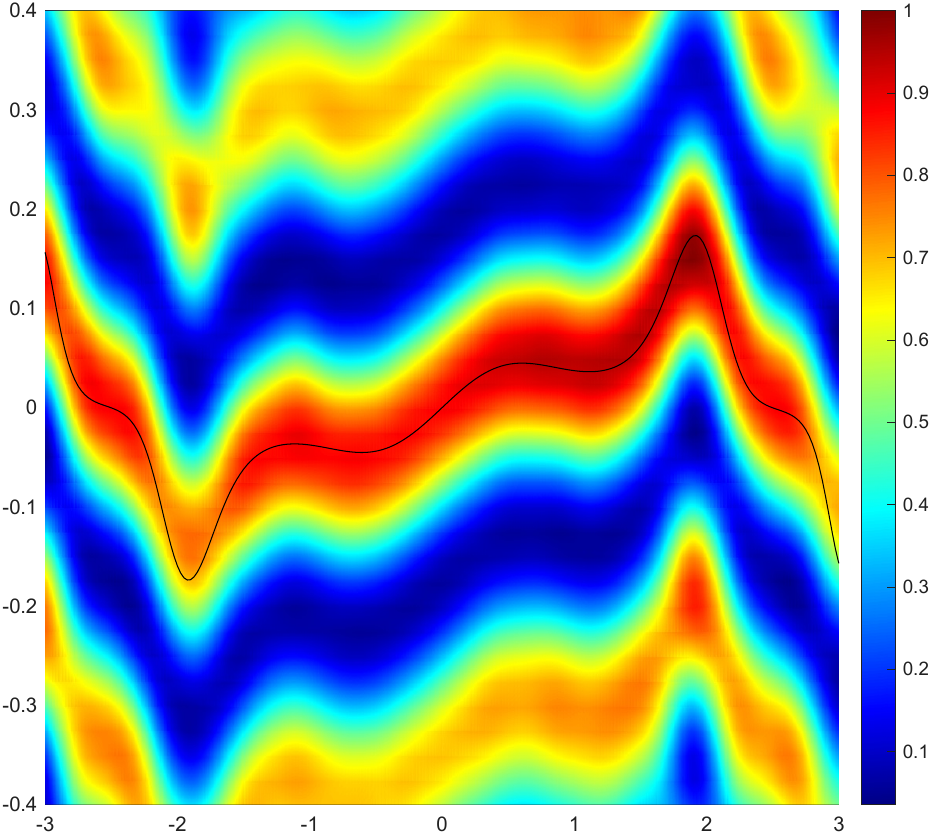}{$\Gamma_6$}
    \caption{TD-DSM imaging results for the six test rough surfaces under the default configuration. The black curves indicate the exact surface profiles.}
    \label{fig:validation}
\end{figure}

Figure~\ref{fig:validation} presents the imaging results, where the exact rough surfaces are superimposed on the corresponding indicator maps. For all six examples, the imaging indicators exhibit strong responses in the vicinity of the true surfaces, successfully identifying their main geometric features, including peaks, valleys, and local variations. Although some oscillatory artifacts appear away from the surfaces, the principal interface structures remain clearly distinguishable. Moreover, all reconstructions are obtained using the same measurement configuration and imaging parameters, without geometry-dependent adjustments.

\setcounter{roughpanel}{0}

These results demonstrate the effectiveness of the proposed TD-DSM in localizing rough surfaces with different geometric characteristics and provide a basis for the subsequent quantitative comparisons.

\subsection{Comparison with frequency-domain DSM}
\label{subsec:frequency-comparison}

\begin{figure}[htbp]
    \centering
    \resultpanelwide{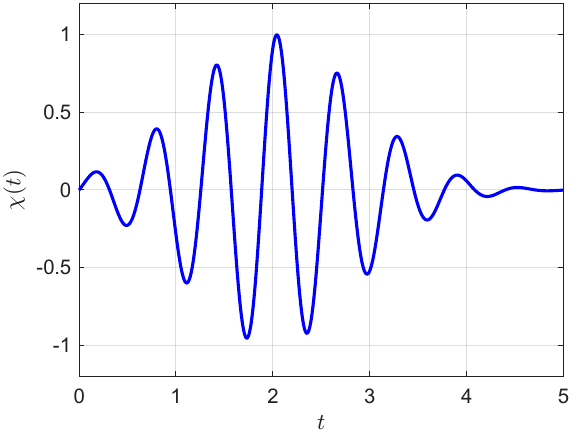}{}\hfill
    \resultpanelwide{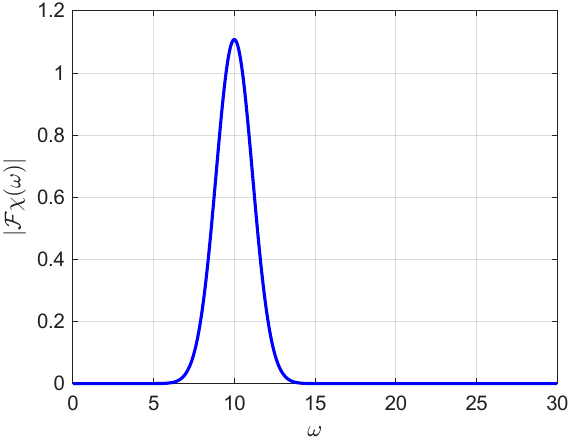}{}
    \caption{The pulse $\chi$ (a) and the magnitude of its Fourier transform $|\mathcal{F}\chi|$ (b).}
    \label{fig:ex_chi}
\end{figure}

To compare the proposed TD-DSM with frequency-domain methods, we retain the baseline pulse with center frequency $k_c=10$. Figure~\ref{fig:ex_chi} shows the pulse $\chi(t;k_c)$ and the magnitude of its Fourier transform. The main spectral energy is concentrated approximately in the frequency range from $8$ to $12$. Thus, we compare TD-DSM with a single-frequency DSM (SF-DSM) at $k=10$ and a multi-frequency DSM (MF-DSM) using $k_m=8,9,\ldots,12$.

The multi-frequency indicator is obtained by summing the single-frequency imaging functions and normalizing the result:
\begin{equation}
    I_{\rm MF}(z)=
    \frac{\sum_{m=1}^{N_k}I_{k_m}(z)}{\max_{z'\in\Omega}\big\{\sum_{m=1}^{N_k}I_{k_m}(z')\big\}}.
    \label{eq:multifrequency-indicator}
\end{equation}
All other measurement parameters remain the same, and the comparison is performed on $\Gamma_1$ and $\Gamma_2$.

\setcounter{widepanel}{0}

\setcounter{roughpanel}{0}
\begin{figure}[htbp]
    \centering
    \resultpanel{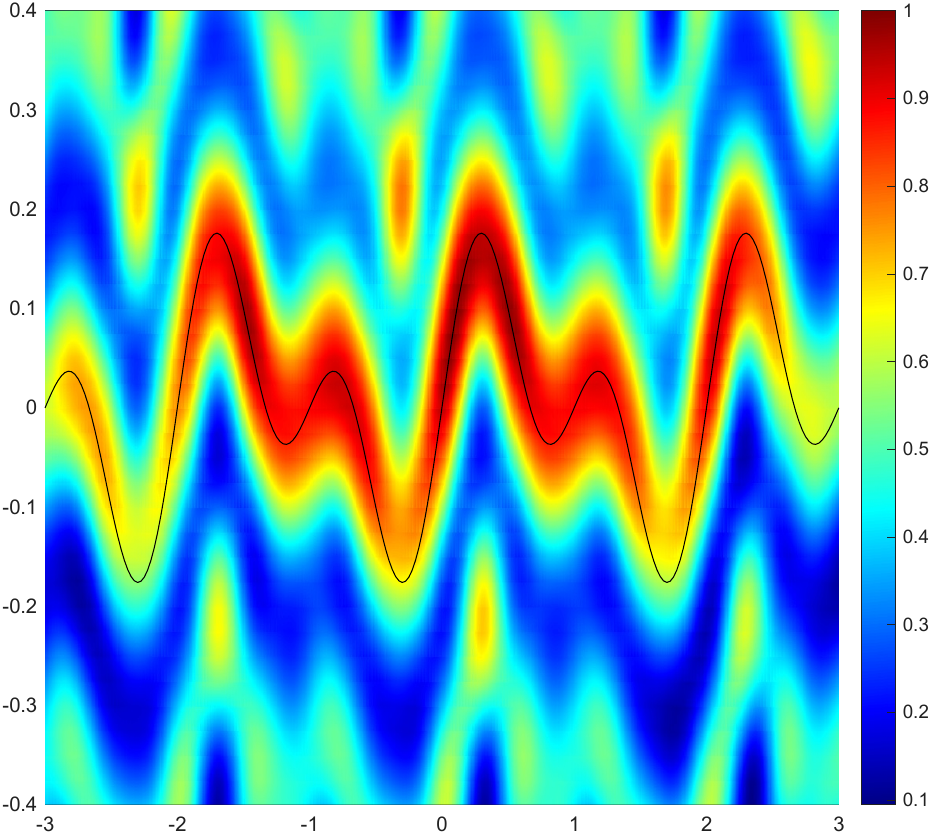}{$\Gamma_1$, SF-DSM}\hfill
    \resultpanel{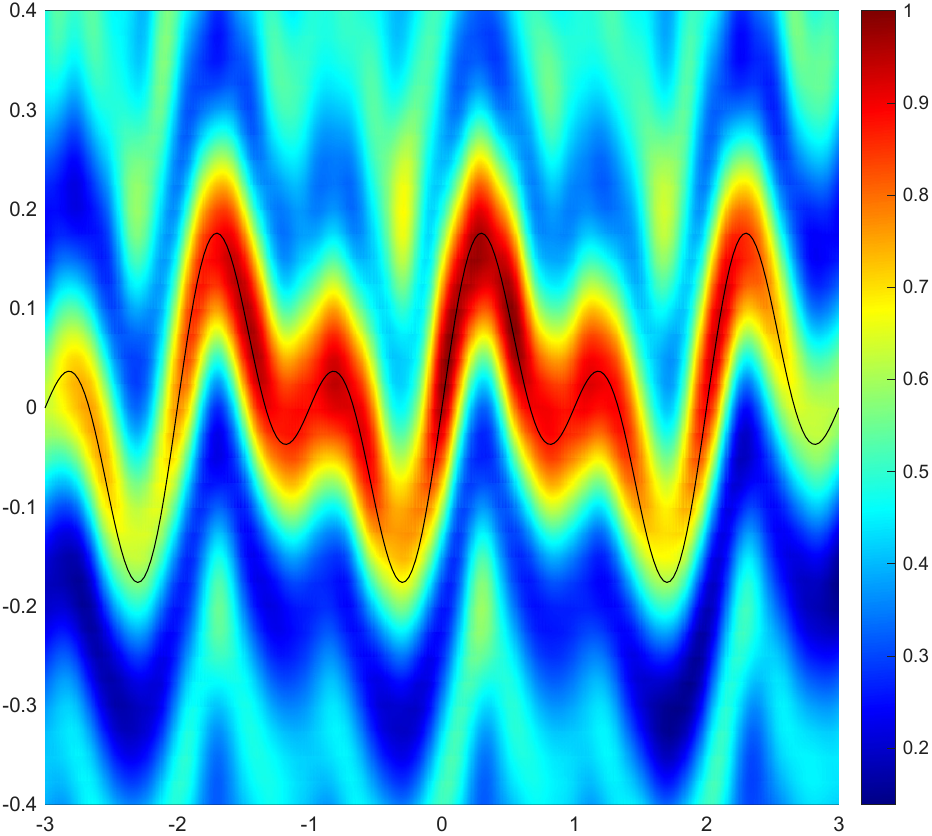}{$\Gamma_1$, MF-DSM}\hfill
    \resultpanel{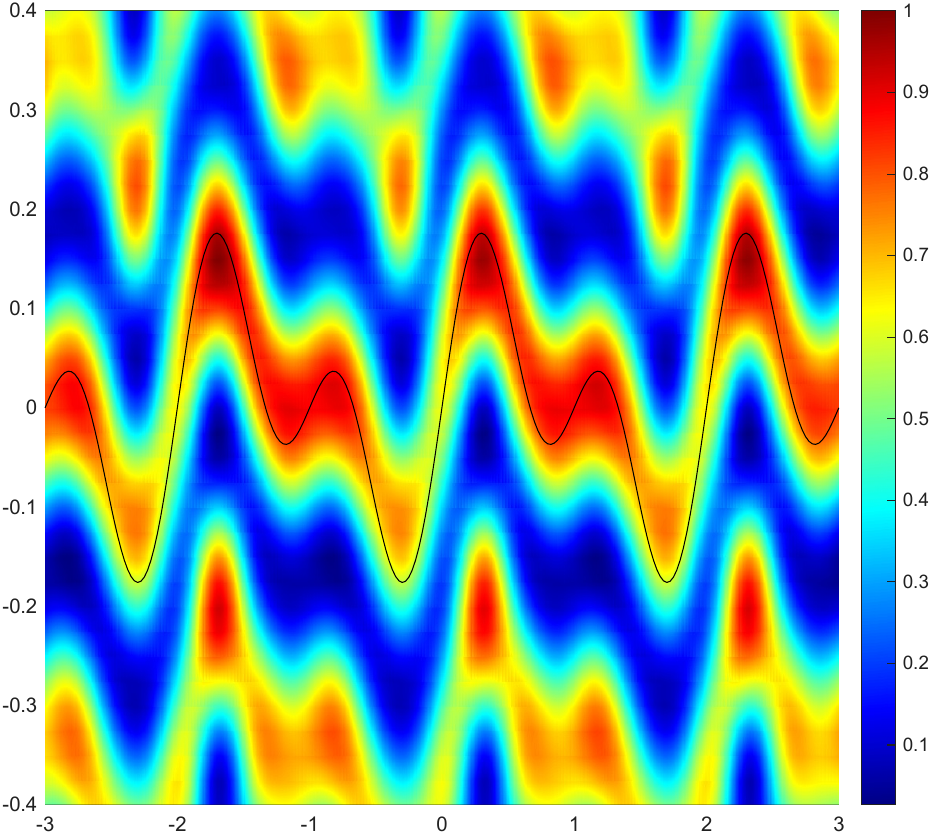}{$\Gamma_1$, TD-DSM}

    \medskip
    \resultpanel{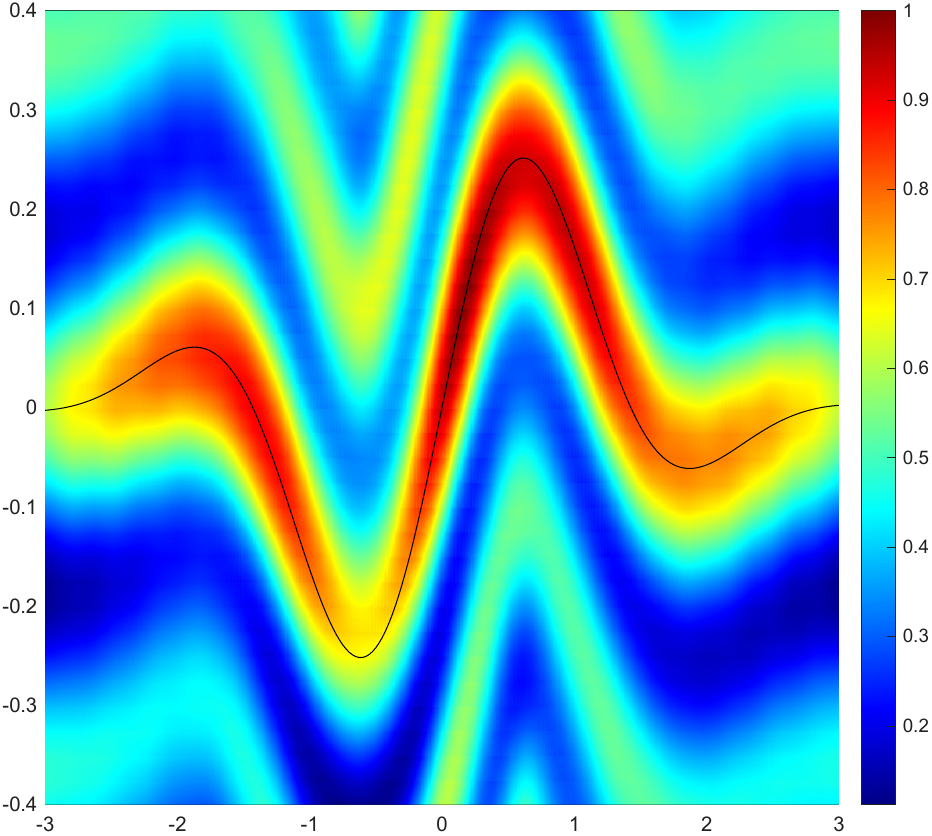}{$\Gamma_2$, SF-DSM}\hfill
    \resultpanel{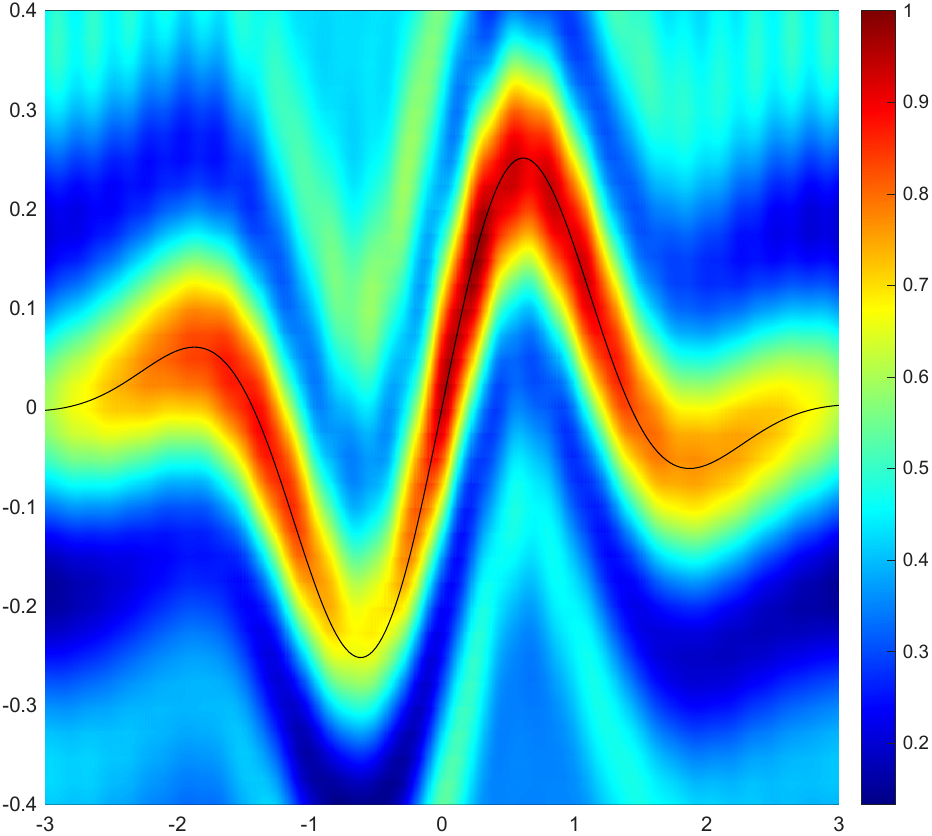}{$\Gamma_2$, MF-DSM}\hfill
    \resultpanel{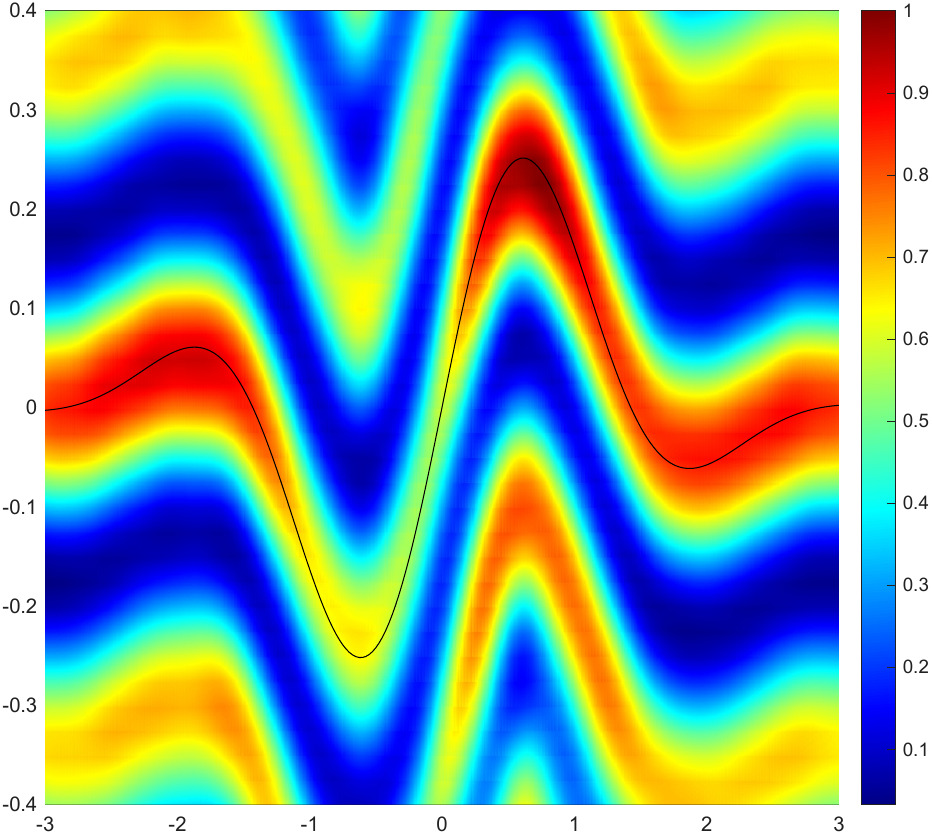}{$\Gamma_2$, TD-DSM}
    \caption{Comparison of SF-DSM, MF-DSM, and TD-DSM for $\Gamma_1$ and $\Gamma_2$.}
    \label{fig:frequency-comparison}
\end{figure}

\begin{table}[htbp]
    \centering
    \caption{Quantitative comparison of SF-DSM, MF-DSM, and TD-DSM for $\Gamma_1$ and $\Gamma_2$.}
    \label{tab:frequency-comparison}
    \begin{tabular}{lcccc}
        \hline
        & \multicolumn{2}{c}{$\Gamma_1$} & \multicolumn{2}{c}{$\Gamma_2$} \\
        \cline{2-3}\cline{4-5}
        Method & $E_{\rm rel}$ & $C$ & $E_{\rm rel}$ & $C$ \\
        \hline
        SF-DSM & 0.9196 & 6.2256 & 0.8809 & 7.2550 \\
        MF-DSM & 0.9180 & 6.2438 & 0.8811 & 7.2790 \\
        TD-DSM & \textbf{0.9131} & \textbf{6.8645} & \textbf{0.8638} & \textbf{8.1902} \\
        \hline
    \end{tabular}
\end{table}

Figure~\ref{fig:frequency-comparison} and Table~\ref{tab:frequency-comparison} show that all three methods can identify the main geometric features of the two rough surfaces. TD-DSM gives a more concentrated response near the true surfaces and better separation from the background, particularly for the localized features of $\Gamma_2$. For $\Gamma_1$, the relative indicator error decreases from $0.9196$ (SF-DSM) and $0.9180$ (MF-DSM) to $0.9131$, while the contrast increases from $6.2256$ and $6.2438$ to $6.8645$. For $\Gamma_2$, the corresponding error decreases to $0.8638$, and the contrast increases to $8.1902$. These results indicate that TD-DSM can effectively use the broadband information of a single pulse and achieve improved localization and imaging contrast in the examples considered here.

\subsection{Resolution analysis}
\label{subsec:resolution}

\subsubsection{Influence of the pulse center frequency}

To examine the effect of pulse center frequency, we reconstruct $\Gamma_3$ using $k_c=5$, $10$, and $15$, with all other parameters fixed at their default values.

\setcounter{roughpanel}{0}
\begin{figure}[htbp]
    \centering
    \resultpanel{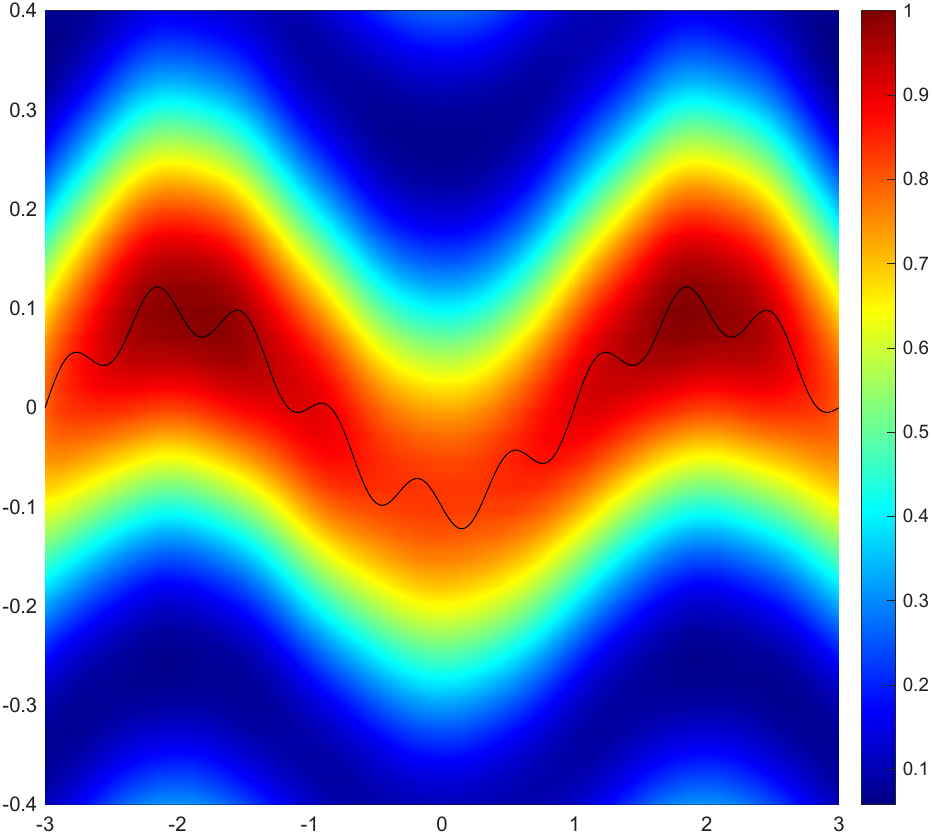}{$k_c=5$}\hfill
    \resultpanel{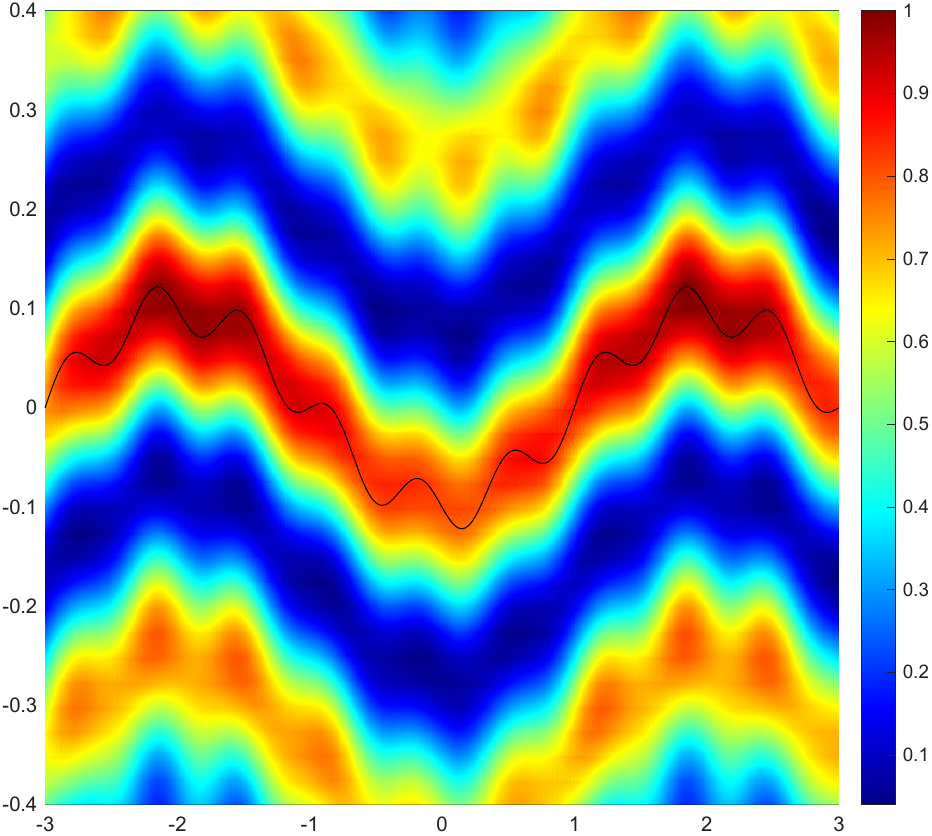}{$k_c=10$}\hfill
    \resultpanel{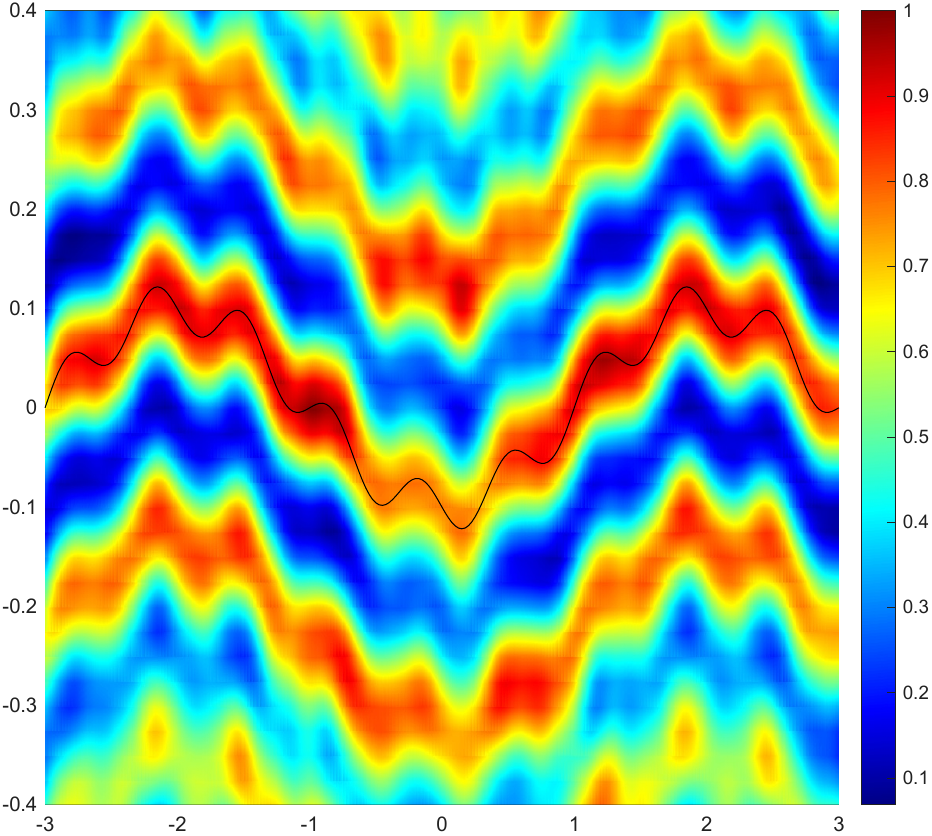}{$k_c=15$}
    \caption{Influence of the pulse frequency on the reconstruction of $\Gamma_3$.}
    \label{fig:pulse-frequency}
\end{figure}

Figure~\ref{fig:pulse-frequency} shows that increasing the pulse center frequency improves the localization of fine-scale features of the rough surface. For $k_c=5$, the indicator response is relatively broad and some small-scale variations are less clearly identified. As $k_c$ increases to $10$ and $15$, the local oscillations of the surface become more distinguishable. However, more pronounced oscillatory responses can also be observed away from the true surface at higher frequencies.

The corresponding values of $E_{\rm rel}$ and $C$ are reported in Table~\ref{tab:pulse-frequency}. As $k_c$ increases from $5$ to $15$, $E_{\rm rel}$ decreases from $0.8871$ to $0.8253$, while $C$ increases from $6.4845$ to $9.8353$. These results indicate that a higher pulse center frequency can improve the localization of fine geometric features in the present configuration, although it does not necessarily suppress all imaging artifacts.

\begin{table}[htbp]
    \centering
    \caption{Influence of the pulse center frequency.}
    \label{tab:pulse-frequency}
    \begin{tabular}{ccc}
        \hline
        $k_c$ & $E_{\rm rel}$ & $C$\\
        \hline
        $5$    & 0.8871 & 6.4845\\
        $10$   & 0.8516 & 8.2952\\
        $15$   & 0.8253 & 9.8353\\
        \hline
    \end{tabular}
\end{table}

\subsubsection{Resolution of closely spaced features}

We next consider a rough surface formed by two Gaussian components,
\begin{equation}
    f_d(x)=0.12e^{-40(x+d/2)^2}+0.12e^{-40(x-d/2)^2}
    \label{eq:two-peak-profile}
\end{equation}
with $d=0.8$, $0.4$, and $0.2$, where $d$ denotes the separation between the centers of the two Gaussian components. The pulse center frequency is $k_c=10$, and all remaining parameters take their default values. This experiment examines the imaging behavior as two localized surface features approach each other.

\setcounter{roughpanel}{0}
\begin{figure}[htbp]
    \centering
    \resultpanel{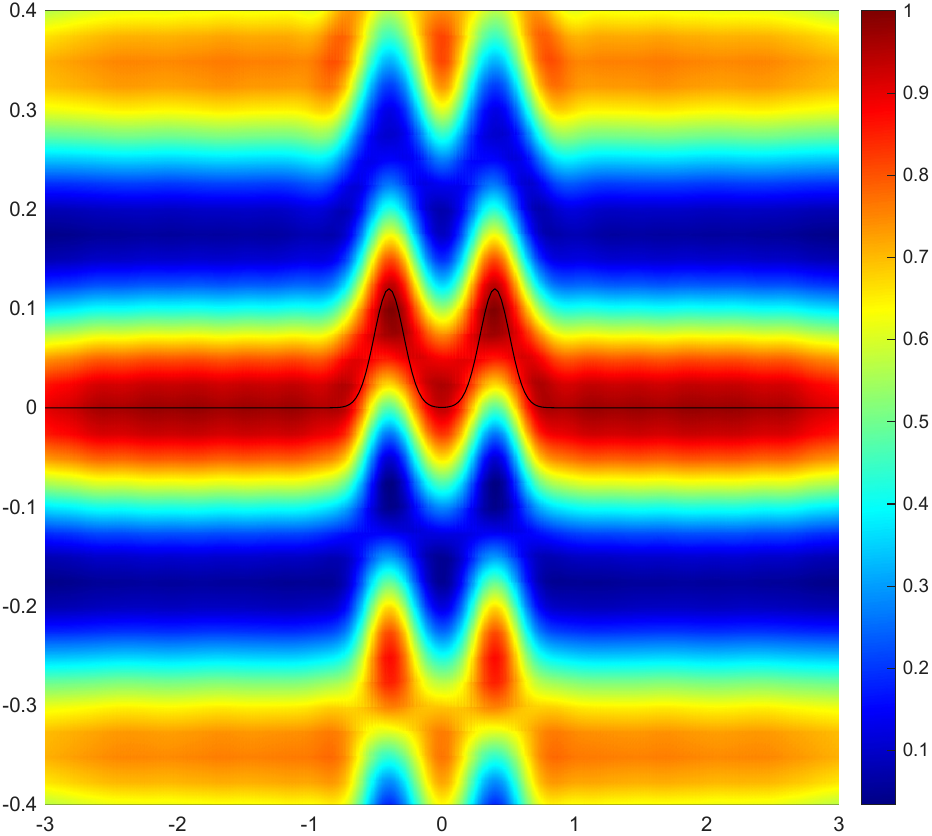}{$d=0.8$}\hfill
    \resultpanel{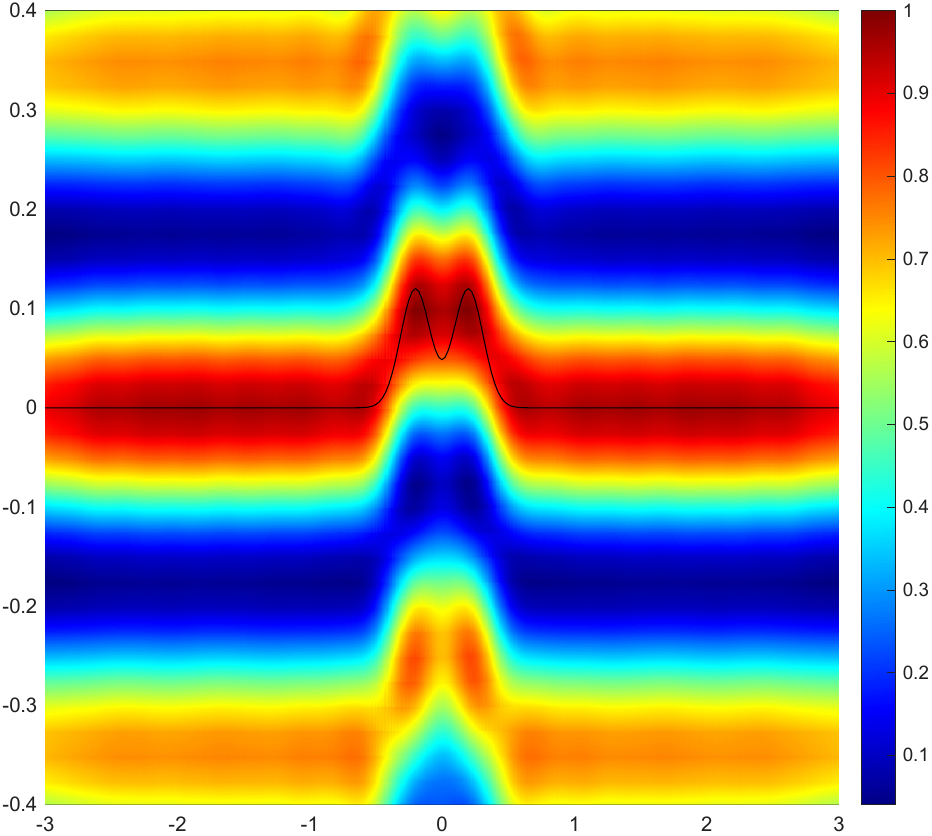}{$d=0.4$}\hfill
    \resultpanel{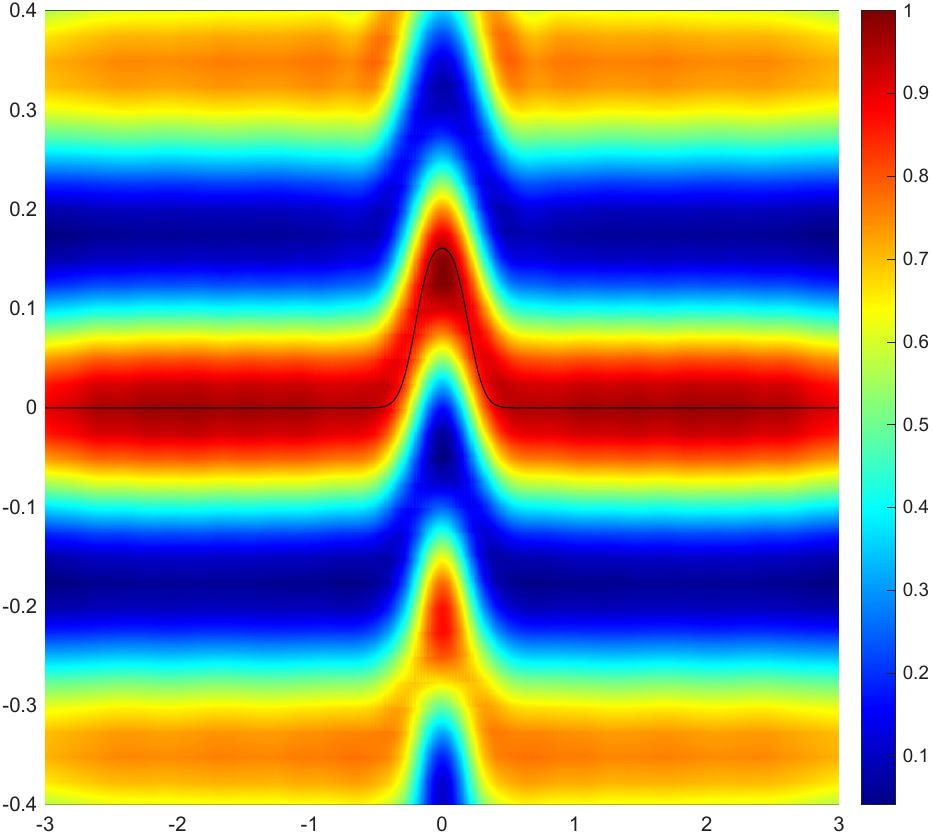}{$d=0.2$}
    \caption{TD-DSM imaging results for two closely spaced Gaussian surface features.}
    \label{fig:two-peak-resolution}
\end{figure}

Figure~\ref{fig:two-peak-resolution} shows that the two localized features are clearly distinguishable for $d=0.8$. When $d=0.4$, the two features are closer and their indicator responses overlap more strongly, making their separation less pronounced. For $d=0.2$, the two Gaussian components merge into a single peak in the exact surface profile, and the imaging result exhibits a corresponding merged response. Thus, the experiment illustrates how the proposed method responds to increasingly close geometric features. In particular, the loss of two distinct peaks at $d=0.2$ should not be interpreted solely as a resolution limitation of the imaging method, since the underlying surface itself has only one local maximum in this case.

\subsection{Influence of measurement configuration}
\label{subsec:measurement-configuration}

\subsubsection{Measurement geometry}

We reconstruct $\Gamma_4$ with apertures $A=3$, $4$, and $5$ at fixed $H=0.8$, and with measurement heights $H=0.6$, $0.8$, and $1$ at fixed $A=4$. All other parameters are kept at their default values.

\setcounter{roughpanel}{0}
\begin{figure}[htbp]
    \centering
    \resultpanel{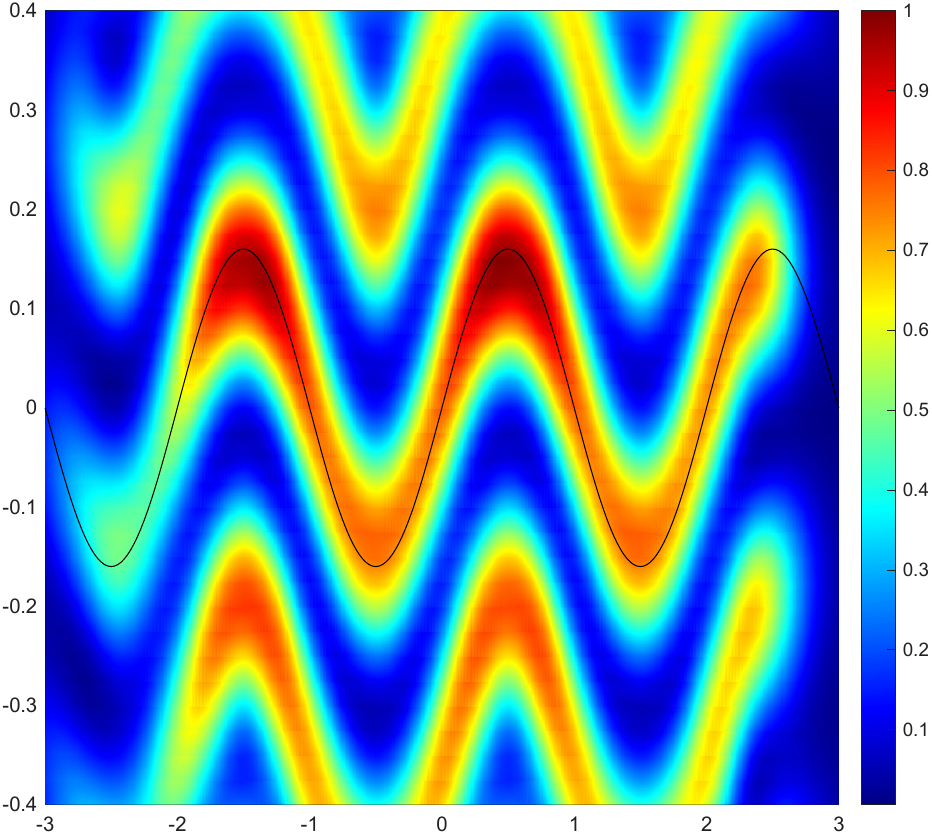}{$A=3$}\hfill
    \resultpanel{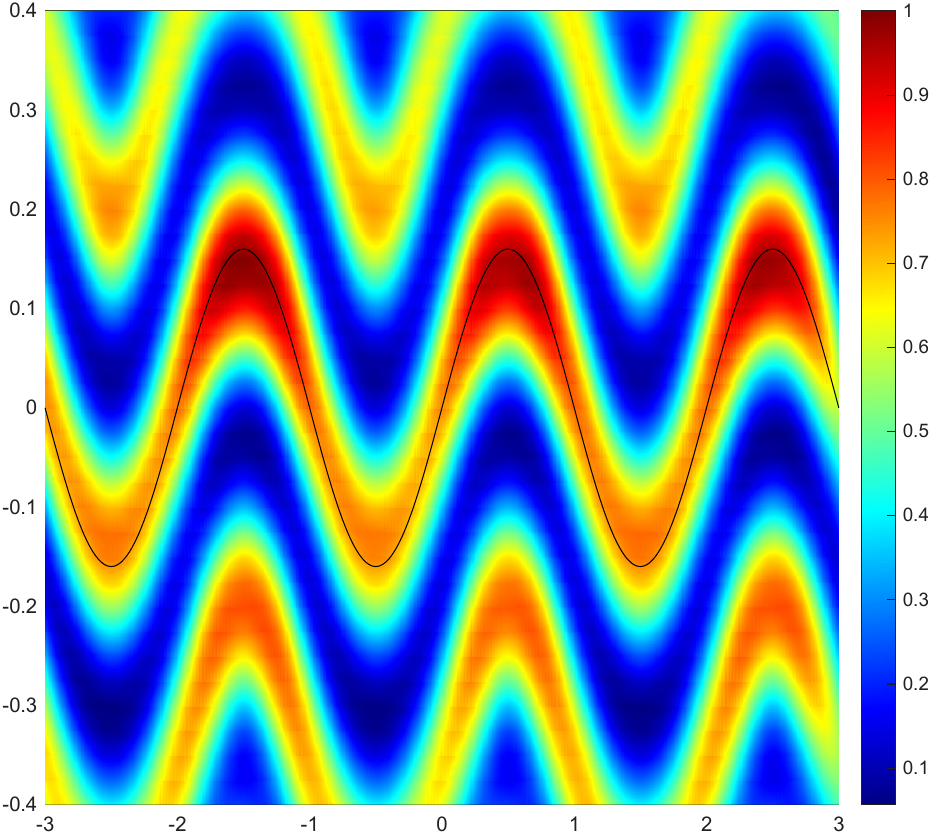}{$A=4$}\hfill
    \resultpanel{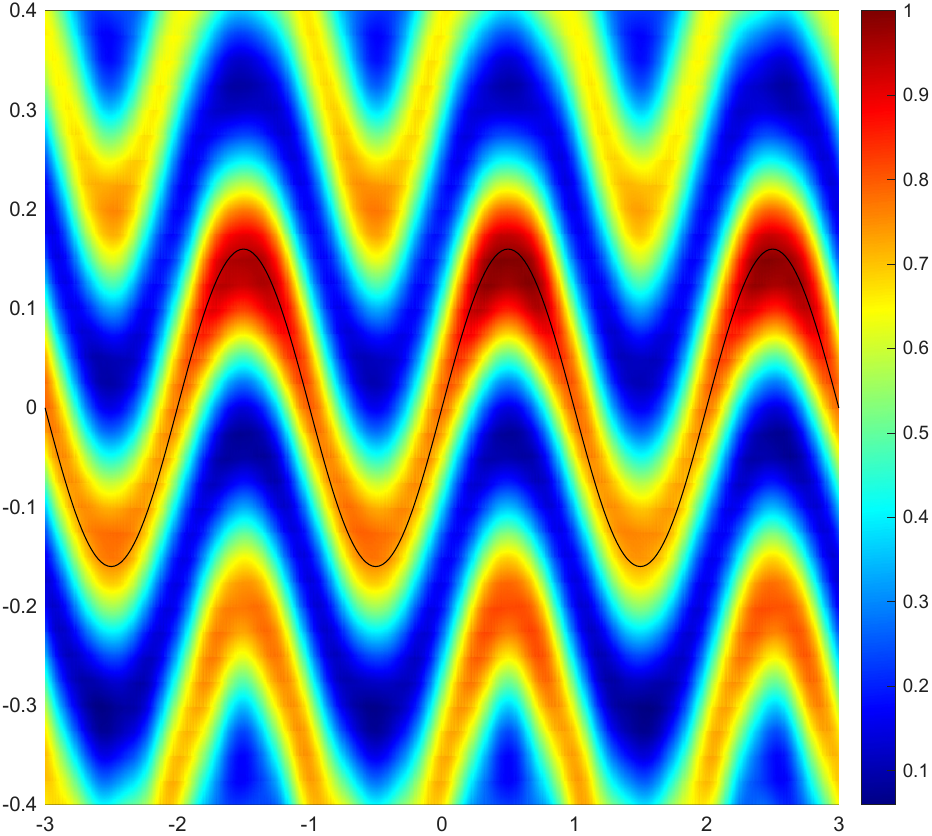}{$A=5$}

    \medskip
    \resultpanel{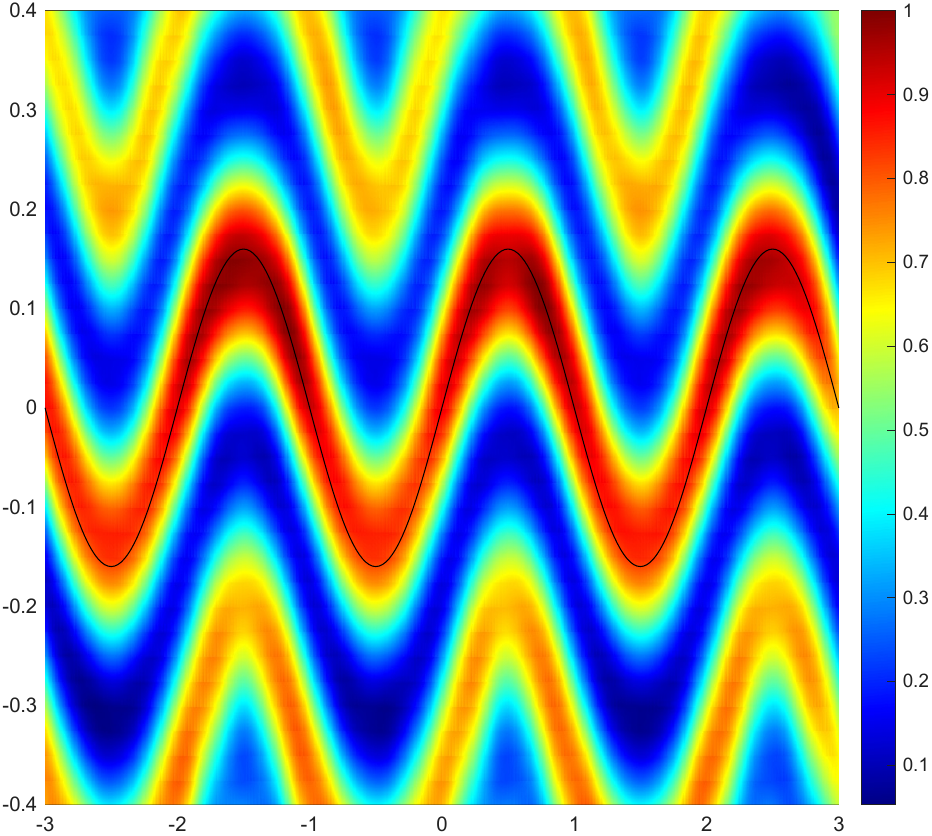}{$H=0.6$}\hfill
    \resultpanel{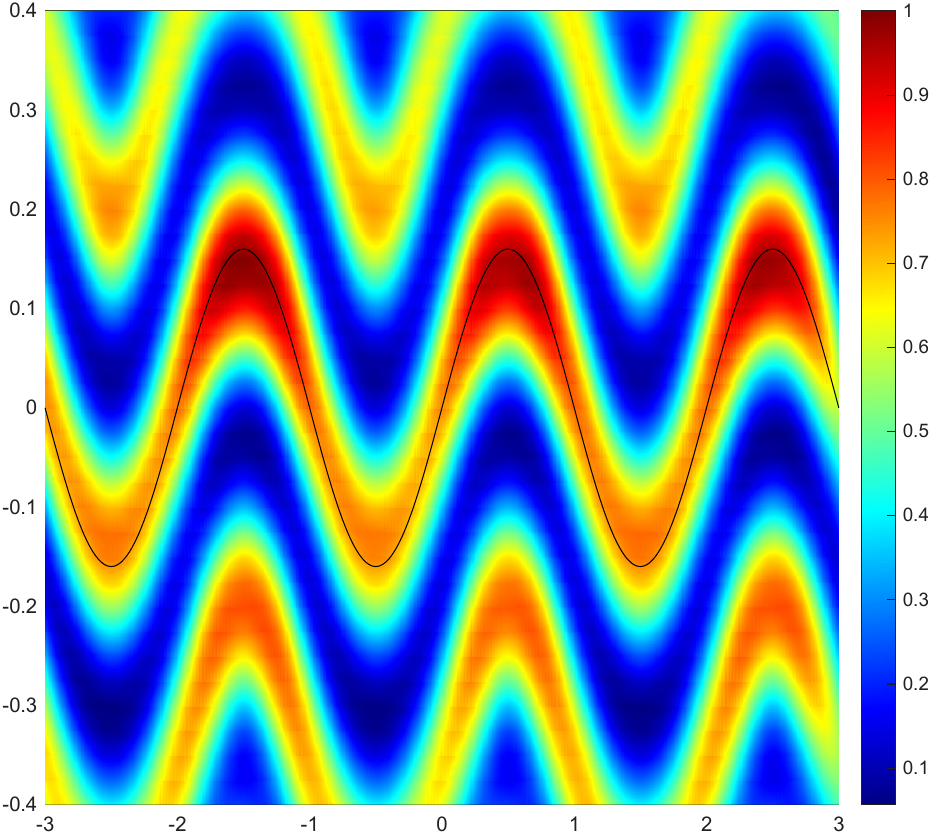}{$H=0.8$}\hfill
    \resultpanel{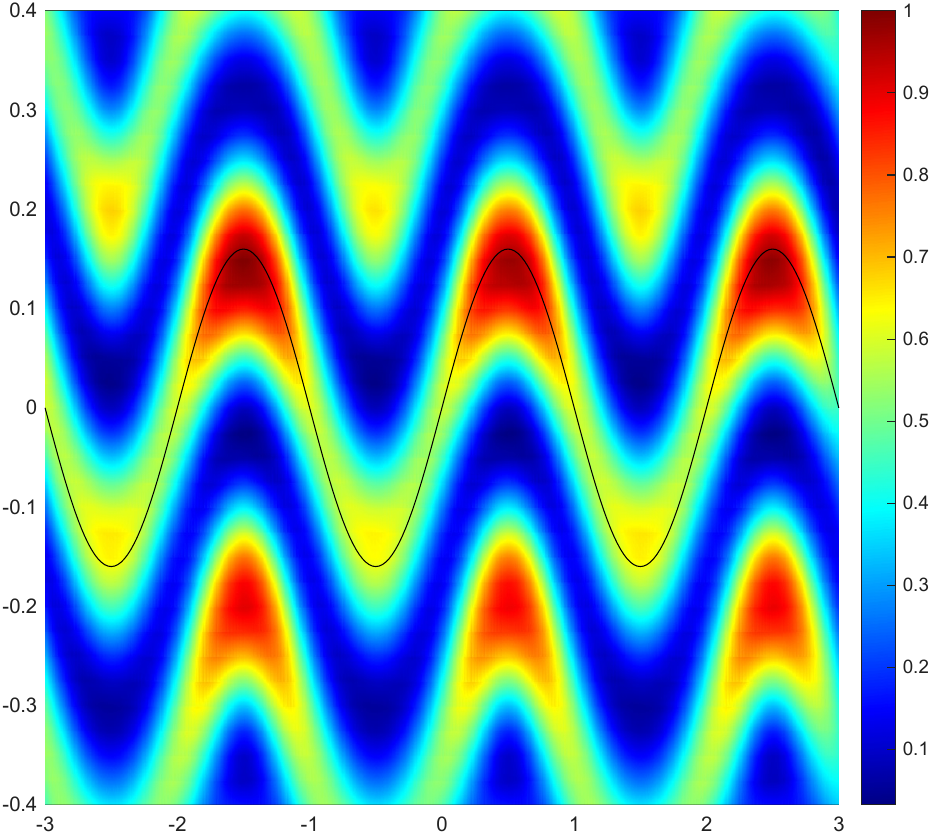}{$H=1$}
    \caption{Influence of aperture parameter at fixed $H=0.8$ (top row) and measurement height at fixed $A=4$ (bottom row) on the reconstruction of $\Gamma_4$.}
    \label{fig:measurement-geometry}
\end{figure}

\begin{table}[htbp]
    \centering
    \caption{Influence of aperture parameter and measurement height.}
    \label{tab:measurement-geometry}
    \begin{tabular}{cccccc}
        \hline
        \multicolumn{3}{c}{Aperture parameter ($H=0.8$)} & \multicolumn{3}{c}{Measurement height ($A=4$)} \\
        \cline{1-3}\cline{4-6}
        $A$ & $E_{\rm rel}$ & $C$ & $H$ & $E_{\rm rel}$ & $C$ \\
        \hline
        $3$ & 0.9613 & 6.7608 & $0.6$ & 0.8734 & 7.4887 \\
        $4$ & 0.8759 & 7.7737 & $0.8$ & 0.8759 & 7.7737 \\
        $5$ & 0.8768 & 7.6825 & $1$ & 0.9017 & 7.9214 \\
        \hline
    \end{tabular}
\end{table}

Figure~\ref{fig:measurement-geometry} shows the influence of the measurement geometry on the imaging results. Increasing the aperture from $A=3$ to $A=4$ improves the localization of the rough surface, whereas a further increase to $A=5$ produces only a limited change. This is also reflected in Table~\ref{tab:measurement-geometry}: $E_{\rm rel}$ decreases from $0.9613$ to $0.8759$ and $C$ increases from $6.7608$ to $7.7737$ when $A$ increases from $3$ to $4$, while both indices remain nearly unchanged for $A=5$.

For the measurement height, a smaller $H$ generally gives a clearer localization of the valley regions. The relative error is lowest at $H=0.6$, while the contrast ratio is highest at $H=1$. Thus, the two indices do not exhibit the same trend, as they measure different aspects of the imaging quality. These results suggest that the measurement geometry affects both the localization and the background response of the indicator.

\subsubsection{Numbers of sources and receivers}

To study the influence of sampling density, we use $\Gamma_5$. In the source study, $N_s=10$, $20$, and $30$ while $N_r=30$; in the receiver study, $N_r=10$, $20$, and $30$ while $N_s=30$. All other parameters remain unchanged.

\setcounter{roughpanel}{0}
\begin{figure}[htbp]
    \centering
    \resultpanel{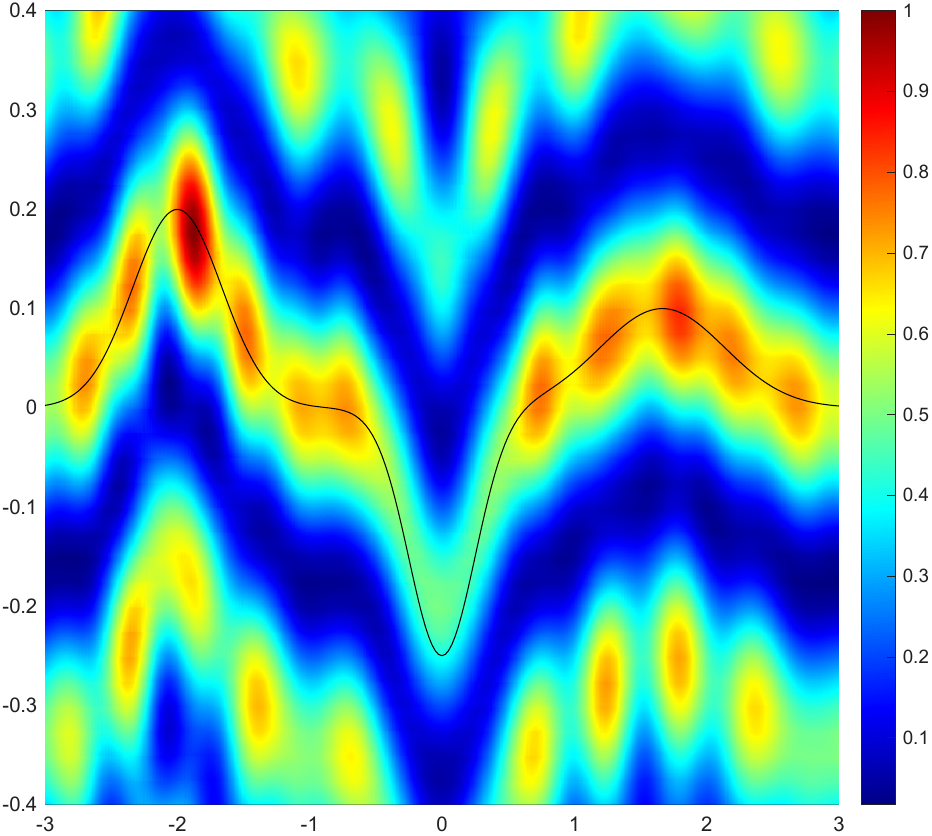}{$N_s=10$}\hfill
    \resultpanel{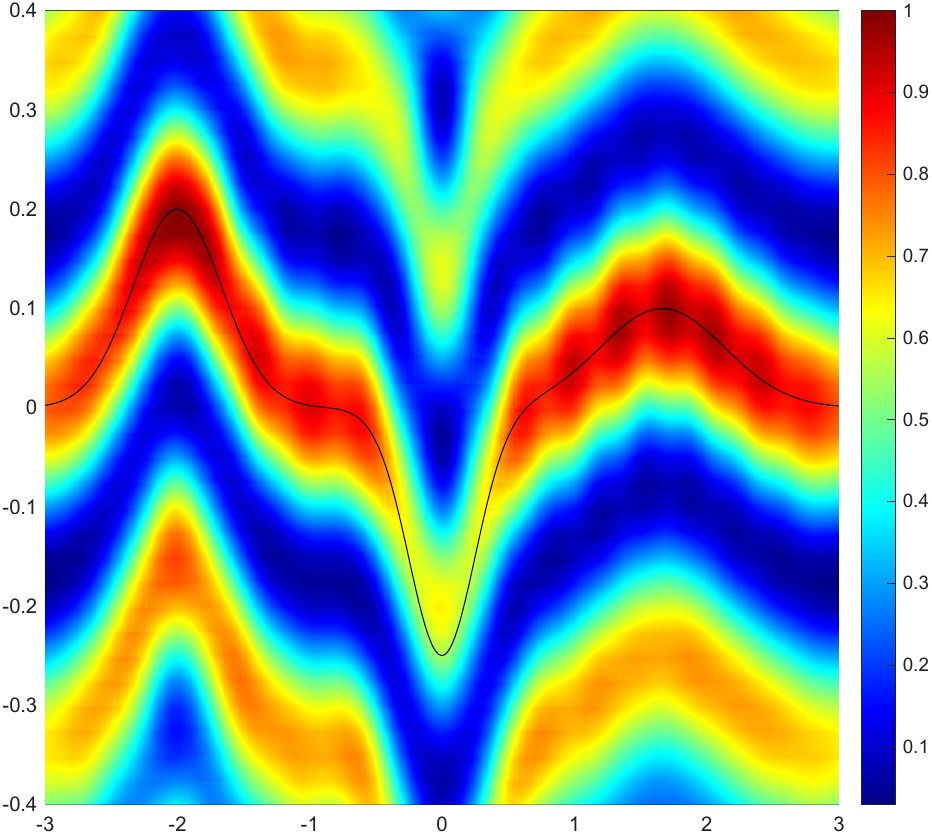}{$N_s=20$}\hfill
    \resultpanel{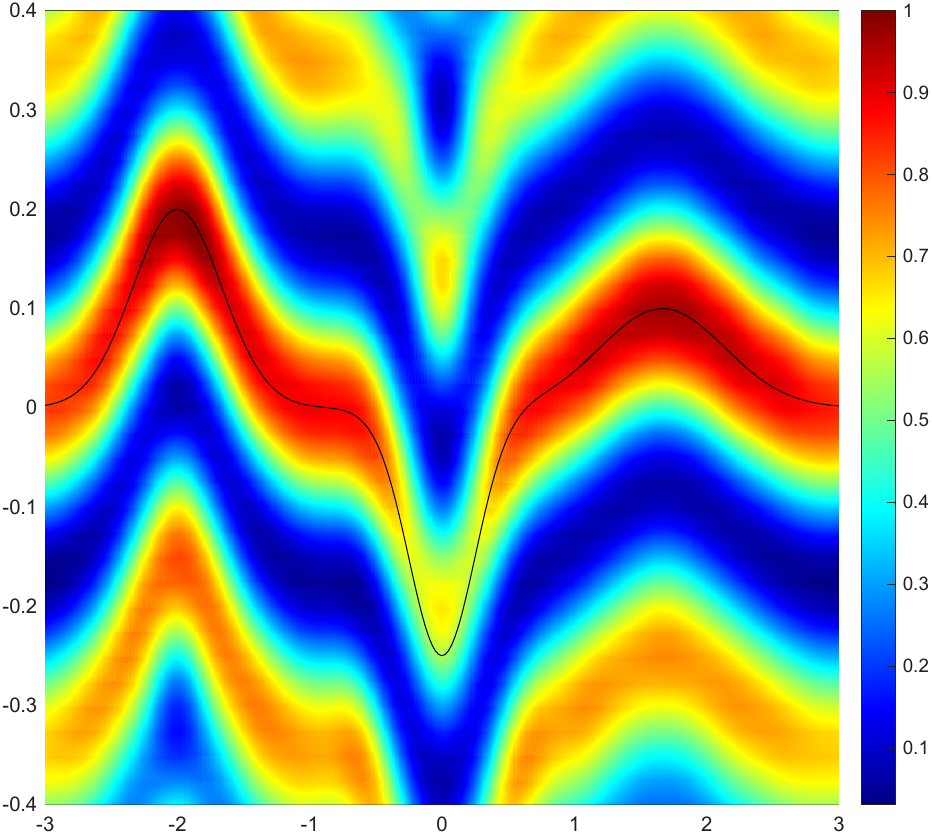}{$N_s=30$}

    \medskip
    \resultpanel{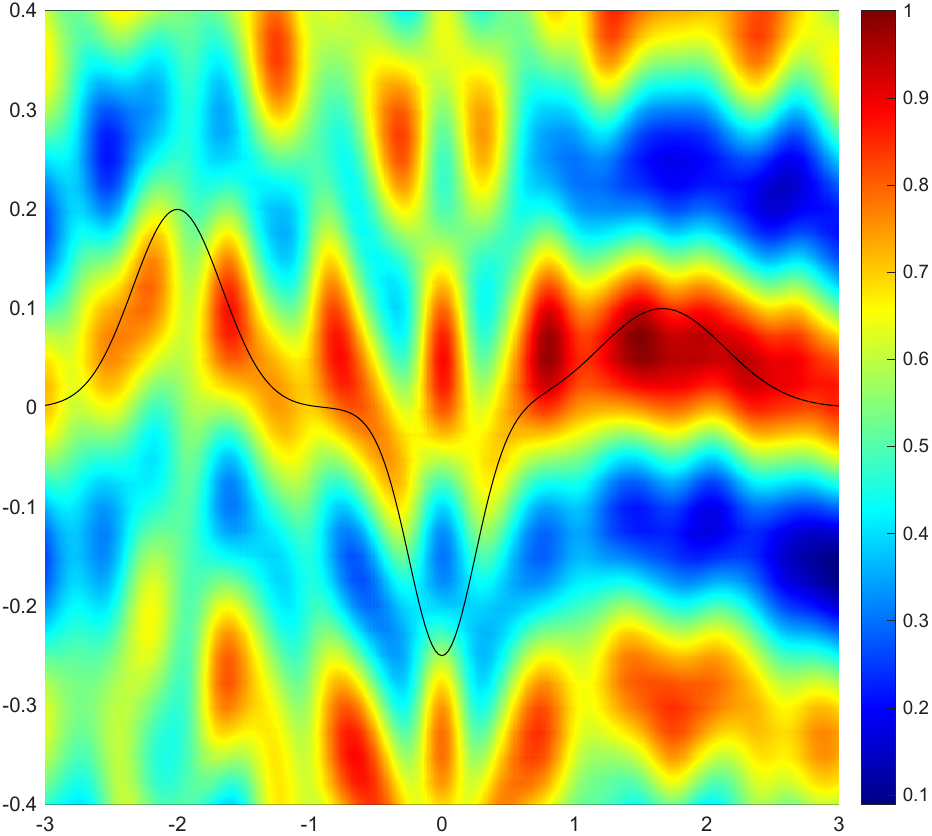}{$N_r=10$}\hfill
    \resultpanel{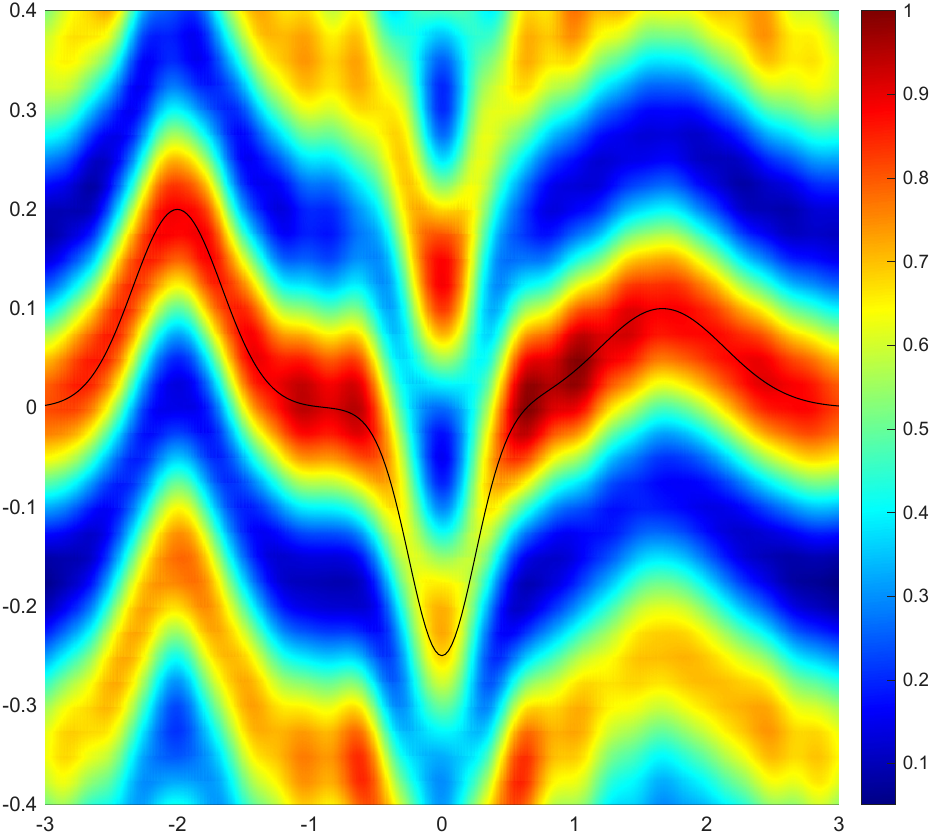}{$N_r=20$}\hfill
    \resultpanel{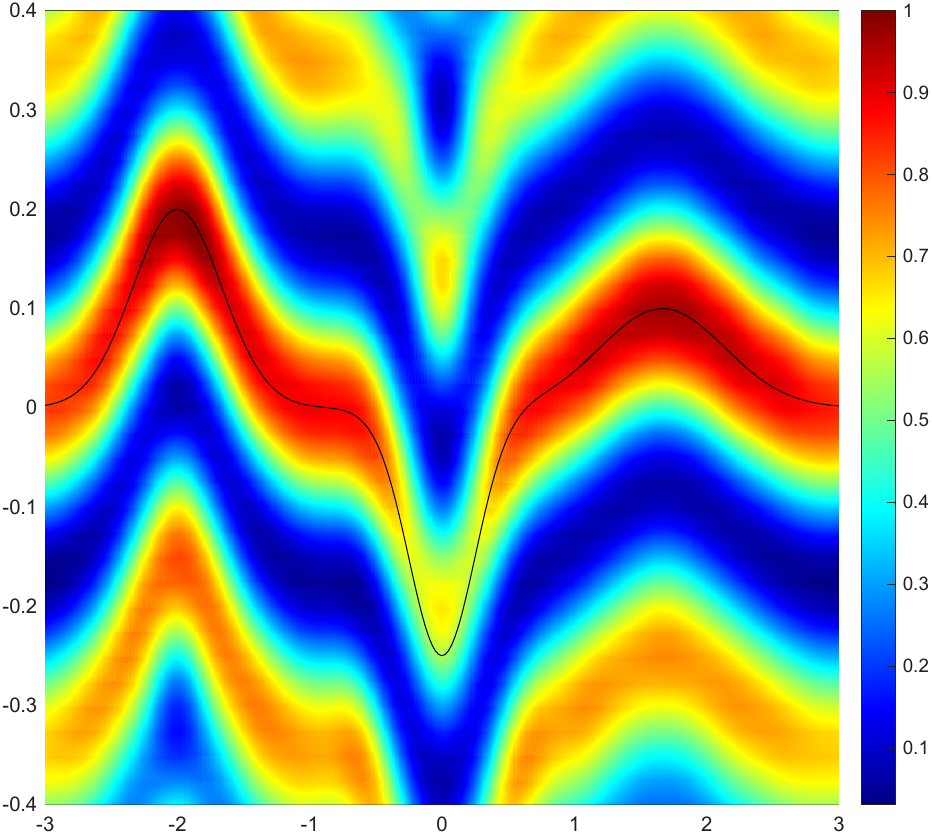}{$N_r=30$}
    \caption{Influence of the number of sources (top row) and receivers (bottom row) on the reconstruction of $\Gamma_5$.}
    \label{fig:sampling-density}
\end{figure}

\begin{table}[htbp]
    \centering
    \caption{Influence of source and receiver sampling density.}
    \label{tab:sampling-density}
    \begin{tabular}{cccccc}
        \hline
        \multicolumn{3}{c}{Number of sources ($N_r=30$)} & \multicolumn{3}{c}{Number of receivers ($N_s=30$)} \\
        \cline{1-3}\cline{4-6}
        $N_s$ & $E_{\rm rel}$ & $C$ & $N_r$ & $E_{\rm rel}$ & $C$ \\
        \hline
        $10$ & 0.9008 & 8.1532 & $10$ & 0.9509 & 5.3998 \\
        $20$ & 0.8600 & 8.2043 & $20$ & 0.8612 & 7.9588 \\
        $30$ & 0.8606 & 8.1603 & $30$ & 0.8606 & 8.1603 \\
        \hline
    \end{tabular}
\end{table}

Figure~\ref{fig:sampling-density} and Table~\ref{tab:sampling-density} show that increasing the number of sources or receivers generally improves the reconstruction when the sampling density is low, but the improvement becomes limited as the number increases. In particular, increasing $N_s$ from $10$ to $20$ reduces $E_{\rm rel}$ from $0.9008$ to $0.8600$, while the contrast ratio remains close to $8.2$. By contrast, increasing $N_r$ from $10$ to $20$ reduces $E_{\rm rel}$ from $0.9509$ to $0.8612$ and increases $C$ from $5.3998$ to $7.9588$. Further increasing either number to $30$ produces only minor changes. In the present configuration, the reconstruction is therefore more sensitive to the receiver sampling density than to the source sampling density.

\subsection{Robustness with respect to noise}
\label{subsec:noise}

Finally, we investigate the robustness of the proposed TD-DSM with respect to measurement noise by reconstructing $\Gamma_6$. For a relative noise level $\delta$, the perturbed scattered-field data are generated by
\begin{equation}
    p^s_\delta=p^s+\delta\,\|p^s\|_\infty X,
    \qquad X\sim U(-1,1),
    \label{eq:noisy-data}
\end{equation}
where $X$ consists of independent samples from the uniform distribution $U(-1,1)$ at the measurement points. Here $\|p^s\|_\infty$ denotes the maximum magnitude of the scattered-field data used in the perturbation. We consider $\delta=0\%$, $20\%$, and $40\%$, while all other parameters remain at their default values.

\setcounter{roughpanel}{0}
\begin{figure}[htbp]
    \centering
    \resultpanel{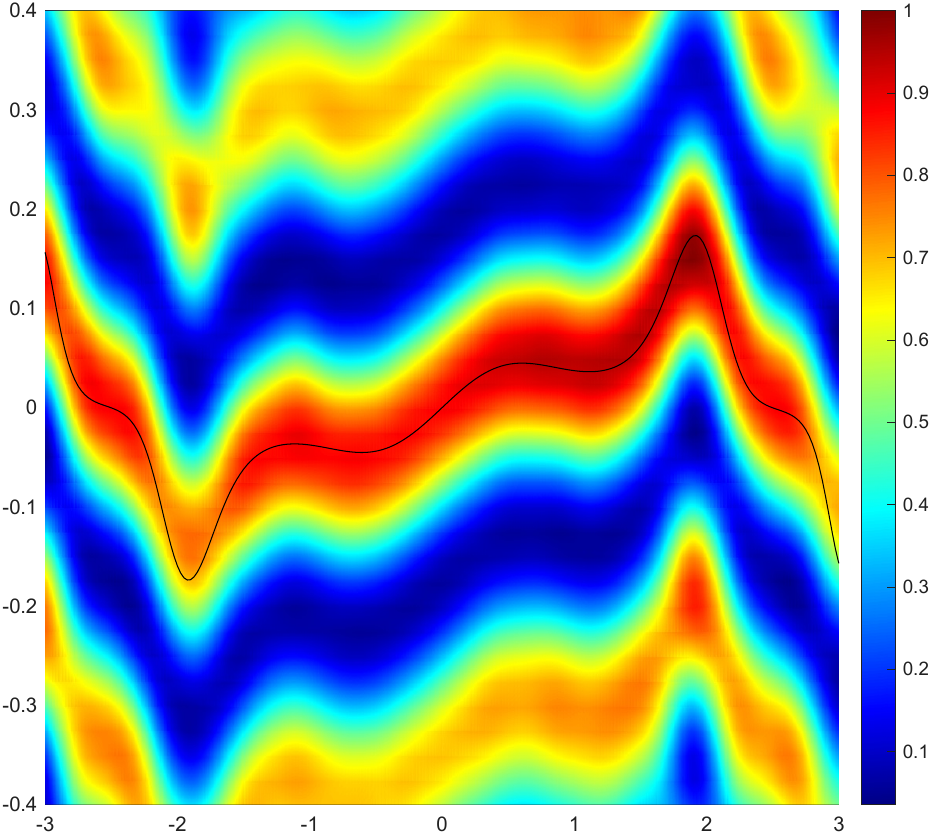}{$0\%$ noise}\hfill
    \resultpanel{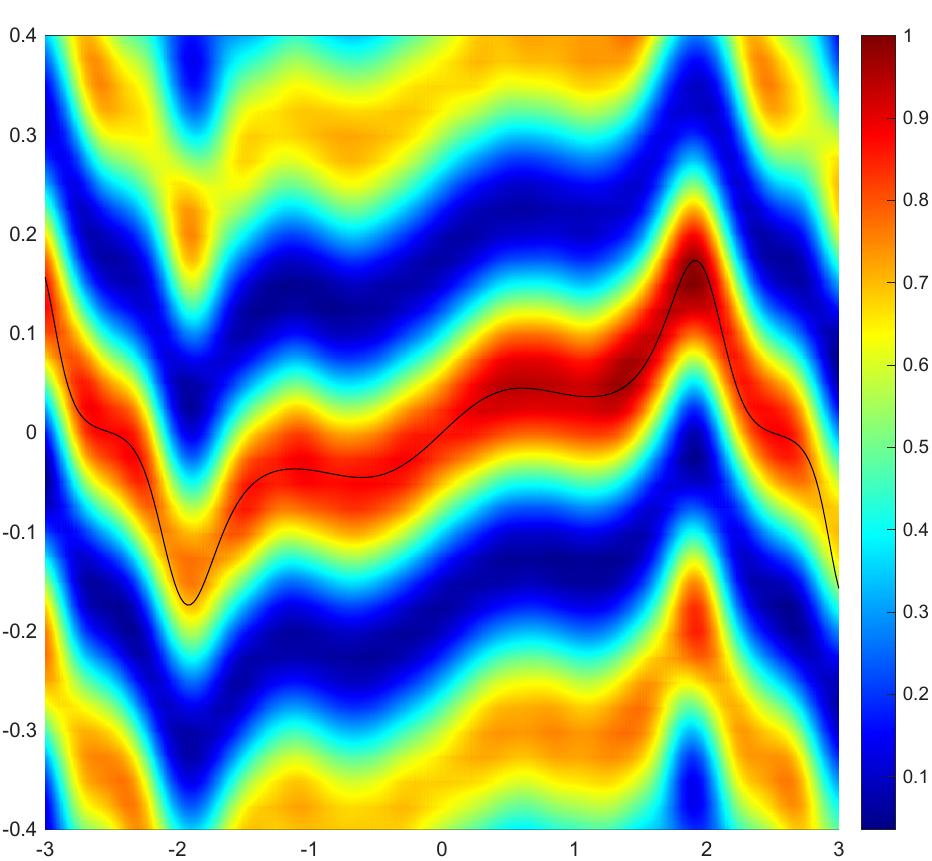}{$20\%$ noise}\hfill
    \resultpanel{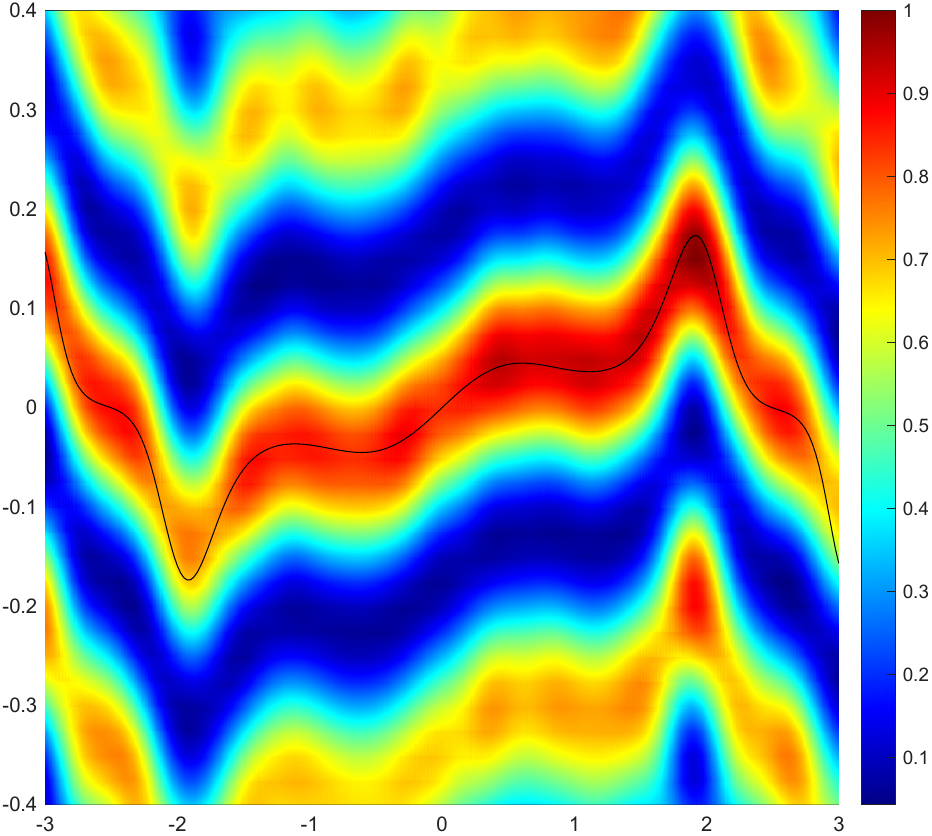}{$40\%$ noise}
    \caption{TD-DSM imaging results for $\Gamma_6$ at different relative noise levels.}
    \label{fig:noise-reconstruction}
\end{figure}

Figure~\ref{fig:noise-reconstruction} shows that the imaging results become less distinct as the noise level increases. Nevertheless, the main geometric features of the rough surface remain identifiable even at $40\%$ noise, although the background responses and local artifacts become more pronounced.

To reduce the influence of random fluctuations, we repeat the noisy experiments with $100$ independent noise realizations and report the averaged values of $E_{\rm rel}$ and $C$ in Table~\ref{tab:noise}.

\begin{table}[htbp]
    \centering
    \caption{Quantitative imaging results for $\Gamma_6$ at different noise levels.}
    \label{tab:noise}
    \begin{tabular}{ccc}
        \hline
        Noise level & $E_{\rm rel}$ & $C$\\
        \hline
        $0\%$  & 0.8603 & 8.0986\\
        $20\%$ & 1.2384 & 6.1840\\
        $40\%$ & 1.4783 & 4.8694\\
        \hline
    \end{tabular}
\end{table}

As shown in Table~\ref{tab:noise}, $E_{\rm rel}$ increases from $0.8603$ to $1.4783$ as the noise level increases from $0\%$ to $40\%$, whereas the contrast ratio decreases from $8.0986$ to $4.8694$. These trends indicate a gradual deterioration in imaging quality under stronger noise contamination. However, the main surface structure can still be localized in the present experiments, demonstrating the robustness of the proposed TD-DSM under the tested noise levels.

\section{Conclusion}

In this paper, we developed a time-domain direct sampling method for reconstructing unbounded rough surfaces from transient near-field measurements. The mathematical analysis provides a theoretical basis for the proposed imaging method, while numerical experiments demonstrate its effectiveness in identifying rough surfaces with different geometric features. The results also suggest several directions for further investigation. In particular, improving the imaging resolution and reducing artifacts under limited-aperture measurements would be of interest, especially for surfaces with fine-scale structures. Another interesting direction is to extend the proposed approach to the inverse scattering problem for periodic surfaces in the time domain. Although periodic surfaces form a special class of rough surfaces, they arise in many important applications. Moreover, their periodic structures may be exploited to develop new imaging methods and further improve the reconstruction performance. These topics will be considered in our future work.



\appendix
\renewcommand{\thetheorem}{\Alph{section}.\arabic{theorem}}
\section{Estimates for a pulsed point source}
\label{app:point_source_estimates}

In this appendix we justify the estimates used in the proof of
Theorem~\ref{thm:well-posedness_reduced}. Let
$\sup f_\Gamma<h<H$, $d=H-h>0$, $s=\omega+\mathrm{i}\sigma$
with $\sigma>0$, and $y=(y_1,H)$. Put
$\breve p^i(x,s)=\breve\chi(s)G_s(x,y)$,
$\beta(\xi,s)=(\xi^2-s^2/c^2)^{1/2}$ with
$\operatorname{Re}\beta>0$, and
$\|v\|_{E_s(D_h)}^2=\|\nabla v\|_{L^2(D_h)^2}^2+
|s|^2c^{-2}\|v\|_{L^2(D_h)}^2$.
All constants below may depend on the fixed positive Laplace weight,
$c$, $d$, and the Lipschitz geometry, but not on $\omega$.

\begin{lemma}\label{lem:point_source_bounds}
For each fixed $\sigma=\sigma_0>0$,
\[
 \|F_s\|_{L^2(\Gamma_h)}
 \le C_{\sigma_0,d,c}(1+|s|)^{1/2}|\breve\chi(s)|,
 \qquad F_s=\partial_\nu \breve p^i-\mathcal B(s)\breve p^i,
\]
and
\[
 \|\breve p^i\|_{E_s(D_h)}
 \le C_{\sigma_0,d,c}(1+|s|)^{1/2}|\breve\chi(s)|.
\]
\end{lemma}

\begin{proof}
With the horizontal Fourier convention of Section 3,
\begin{equation*}
 \widehat{\breve p^i}(\xi,x_2,s)=
  \frac{\breve\chi(s)}{2\beta(\xi,s)}
  e^{-\beta(\xi,s)(H-x_2)}e^{\mathrm{i}\xi y_1},
  \qquad x_2<H.
\end{equation*}
Hence $\widehat F_s(\xi)=\breve\chi(s)
e^{-d\beta(\xi,s)}e^{\mathrm{i}\xi y_1}$.
Set $k=|\omega|/c$, $a=\sigma_0/c$,
$b=\operatorname{Re}\beta>0$, and $R=|\beta|^2$.
The identities
\[
 R^2=(\xi^2-k^2+a^2)^2+4k^2a^2,\qquad
 b^2=\frac{R+\xi^2-k^2+a^2}{2},\qquad b\ge a
\]
imply, for $|\xi|>2k+1$, that
$b\ge\sqrt3|\xi|/2$.
Consequently,
\[
 \int_{\mathbb R}e^{-2db}\,\mathrm d\xi
 \le (4k+2)e^{-2da}
 +2\int_{2k+1}^\infty e^{-\sqrt3d\xi}\,\mathrm d\xi
 \le C_{a,d}(1+k).
\]
Parseval's identity proves the first estimate.

Since $D_h\subset\{x_2<h\}$, integrating the squared
Fourier representation over the lower half-plane gives
\[
 \|\breve p^i\|_{E_s(D_h)}^2
 \le\frac{|\breve\chi(s)|^2}{16\pi}\,
 \mathcal I(k),\qquad
 \mathcal I(k)=\int_{\mathbb R}
 \frac{\xi^2+R+k^2+a^2}{Rb}e^{-2db}\,\mathrm d\xi.
\]
Using $b^2-(\operatorname{Im}\beta)^2
=\xi^2-k^2+a^2$, we have the exact identity
\[
 \xi^2+R+k^2+a^2=2(b^2+k^2).
\]
We show $\mathcal I(k)\le C_{a,d}(1+k)$.
For $k\ge\max\{1,2a\}$ and $0\le\xi\le2k$,
\[
 R\ge Ck(|\xi-k|+a),\qquad
 \frac bR\le\frac1a,\qquad
 \frac{k^2}{Rb}\le\frac{k}{a\sqrt R}
 \le C_a\frac{\sqrt k}{\sqrt{|\xi-k|+a}}.
\]
Here the middle inequality follows from
$|\operatorname{Re}\beta\,\operatorname{Im}\beta|=ka$.
It follows that
\[
 \int_0^{2k}\frac{2(b^2+k^2)}{Rb}
 e^{-2db}\,\mathrm d\xi
 \le C_a\left(k+\sqrt k
 \int_0^{2k}\frac{\mathrm d\xi}
 {\sqrt{|\xi-k|+a}}\right)
 \le C_a k.
\]
For $\xi>2k$, $b\ge\sqrt3\xi/2$ and
$R\ge3\xi^2/4$, so, since $k\ge1$,
\[
 \int_{2k}^\infty
 \frac{2(b^2+k^2)}{Rb}e^{-2db}\,\mathrm d\xi
 \le C_a\int_{2k}^\infty
 (\xi^{-1}+\xi^{-3})e^{-\sqrt3d\xi}\,\mathrm d\xi
 \le C_{a,d}.
\]
The negative half-line is identical. When $k$ is
bounded, $b\ge a$ and the exponential decay for
large $|\xi|$ give a uniform bound directly.
Thus $\mathcal I(k)\le C_{a,d}(1+k)$, proving
the second estimate.
\end{proof}

\paragraph{Causality and inversion}
We verify the analyticity and polynomial growth conditions required for the inverse Laplace transform. For $\mathrm{Im}\,s>0$, the DtN operator $\mathcal B(s)$ and the point-source data depend analytically on $s$. Since the associated variational operator is coercive and depends analytically on $s$, the analytic inverse operator theorem implies that $\breve p^s$ is analytic in the upper half-plane.

It remains to establish a polynomial growth estimate uniformly for $\mathrm{Im}\,s\geq\sigma_0>0$. Using the Fourier representation in the proof of Lemma~\ref{lem:point_source_bounds}, together with $\mathrm{Re}\,\beta(\xi,s)\geq\sigma/c$, and splitting the Fourier integrals into $|\xi|\leq 2|s|/c+1$ and its complement, we obtain
\[
\begin{aligned}
\|F_s\|_{L^2(\Gamma_h)}
&\leq C_{\sigma_0}(1+|s|)^{1/2}|\breve\chi(s)|,\\
\|\breve p^i\|_{E_s(D_h)}
&\leq C_{\sigma_0}(1+|s|)^{3/2}|\breve\chi(s)|,
\end{aligned}
\]
uniformly for $\mathrm{Im}\,s\geq\sigma_0$. Indeed, the lower bound on $\mathrm{Re}\,\beta$ controls the integrands in the first region, while the exponential decay for large $|\xi|$ controls the remaining integrals.

The variational estimate used in the proof of Theorem~\ref{thm:well-posedness_reduced} remains valid for $\sigma\geq\sigma_0$. Combining this estimate with the above bounds and $\breve p^s=\breve p-\breve p^i$, we obtain
\[
\|\breve p^s(\cdot,s)\|_{H^1(D_h)}
\leq C_{\sigma_0}(1+|s|)^2|\breve\chi(s)|.
\]
Since $\chi$ is causal, the Cauchy--Schwarz inequality gives
\[
|\breve\chi(s)|\leq (2\sigma_0)^{-1/2}\|\chi\|_{L^2(\mathbb R)},\qquad \mathrm{Im}\,s\geq\sigma_0.
\]
Consequently,
\[
\|\breve p^s(\cdot,s)\|_{H^1(D_h)}
\leq C_{\sigma_0}(1+|s|)^2\|\chi\|_{H^1(\mathbb R)},\qquad \mathrm{Im}\,s\geq\sigma_0.
\]

Thus, the analyticity and polynomial growth conditions for the vector-valued inverse Laplace transform are satisfied. By the corresponding Paley--Wiener--Laplace theorem, $\breve p^s$ is the Laplace transform of a causal $H^1(D_h)$-valued distribution.

\end{document}